\documentclass{amsart}
\usepackage{graphicx}

\newtheorem{thm}{Theorem}[section]

\theoremstyle{definition}

\theoremstyle{remark}

\numberwithin{equation}{section}

\begin{document}

\title[Numerical study of moments of the zeta function]{The zeta function on the critical line: \\
numerical evidence for moments and random matrix theory models}

\author[G.A. Hiary]{Ghaith A. Hiary}
\thanks{Preparation of this material was partially supported by the National Science Foundation 
under agreements No.  DMS-0757627 (FRG grant) and DMS-0635607.  Computations were carried out
at the Minnesota Supercomputing Institute.}
\address{Pure Mathematics, University of Waterloo, 200 University Ave West, Waterloo, Ontario, Canada, N2L 3G1.}
\email{hiaryg@gmail.com}

\author[A.M. Odlyzko]{Andrew M. Odlyzko}
\thanks{}
\address{Department of Mathematics, University of Minnesota, 206 Church St. S.E., Minneapolis, MN, 55455.}
\curraddr{}
\email{odlyzko@umn.edu}

\subjclass[2000]{Primary, Secondary}


\begin{abstract}
Results of extensive computations of moments of the Riemann 
zeta function on the critical line are presented. 
Calculated values are compared with predictions motivated by random 
matrix theory.  The results can help in deciding
between those and competing predictions. It is shown that for high moments and at large heights, the variability
of moment values over adjacent intervals is substantial,
even when those intervals are long, as long as a block
containing $10^9$ zeros near zero number $10^{23}$.
More than
anything else, the variability
illustrates the limits of what one can learn about the zeta function from numerical
evidence. 

It is shown the rate of decline of extreme values of the 
moments is modeled relatively well by power laws. 
Also, some long range correlations
in the values of the second moment, as well as asymptotic oscillations
in the values of the shifted fourth moment, are found.

The computations described here relied on several
representations of the zeta function.  The numerical comparison 
of their effectiveness that is presented is of independent interest,
for future large scale computations.
\end{abstract}

\maketitle

\section{Introduction}

Absolute moments of the Riemann zeta function on the critical line have been 
the subject of intense theoretical investigations by 
Hardy, Littlewood, Selberg, Titchmarsh, and many others. 
It has long been  conjectured 
that the $2k^{th}$ moment of $|\zeta(1/2+it)|$ should grow like $c_k (\log t)^{k^2}$ for 
some constant $c_k >0$. Conrey and Ghosh~\cite{CG1} reformulated this long-standing conjecture 
in a precise form: for a fixed integer $k>0$, 

\begin{equation}  \label{eq:conjmoments}
\frac{1}{T}\int_0^T |\zeta(1/2+it)|^{2k}\, dt \sim \frac{a(k)g(k)}{k^2!}\left(\log T\right)^{k^2}\, ~~as~~~~T \to \infty,
\end{equation}

\noindent
where $a(k)$ is a certain, generally understood, ``arithmetic factor,'' and $g(k)$ is an integer.
 Keating and Snaith~\cite{KS} suggested 
 there might be relations between random matrix theory (RMT) and moments of the zeta function, which led them to a precise conjecture for $g(k)$. 

More generally, it is expected the $2k^{th}$ moment of $|\zeta(1/2+it)|$ should grow like

\begin{equation} \label{eq:rubmomentold10}
\frac{1}{T}\int_0^T |\zeta(1/2+it)|^{2k}\,dt\sim \frac{1}{T}\int_0^T P_k\left(\log \frac{t}{2\pi}\right) \,dt\,,
\end{equation}

\noindent
where $P_k(x)$ is a polynomial of degree $k^2$ with leading coefficient $a(k)g(k)/k^2!$, which is the same leading coefficient as in (\ref{eq:conjmoments}). 
 Recently, Conrey, Farmer, Keating, Rubinstein, and
Snaith~\cite{CFKRS1} gave a conjectural recipe, compatible with the Keating-Snaith prediction for $g(k)$,
to compute lower order coefficients of $P_k(x)$;\footnote{See also Diaconu, Goldfeld, and Hoffstein~\cite{DGH}.} see Section 2.

The purpose of this article is to present the results of  
numerical computations of integral moments of $|\zeta(1/2+it)|$ at relatively large heights.
The main goal is to enable comparison with 
RMT, and other more complete predictions for moments of the zeta function on the critical line. These predictions are used to get insights that go beyond
what can be derived even with the Riemann hypothesis. 

RMT has produced a variety of (so far still uproven) conjectures
for the zeta function.  Many of these conjectures have been
either inspired by, or 
reinforced by, numerical evidence.  In general, the agreement
between numerical data and RMT predictions has been very good.
But it is also known there are some quantities,
such as the number variance in the distribution of spacings between zeta zeros, that depend explicitly on
primes.  Therefore, for some ranges of the spacing parameter,
these quantities do not follow RMT predictions exactly.

Similarly, there are suggestions that
 the behavior of high moments of the zeta function may involve
arithmetical factors such as $a(k)$ in conjecture (\ref{eq:conjmoments}), and so might not follow the RMT predictions completely.
We provide numerical
data that might be used in the future to shed light on this
issue.  Even today, the statistics we present on the
variation in the moment values over various intervals might be informative
in judging the extent to which RMT provides an adequate explanation
for observed behavior.

We computed the moments of $|\zeta(1/2+it)|$  
for a set of 15 billion zeros near $t=10^{22}$, and for sets of one billion zeros each near $t=10^{19}$ and $t=10^{15}$. We also computed the moments near $t=10^7$.
These zero samples were used because they were the main ones
available from prior computations by the second author. 

A relatively large sample size (such as a billion zeros) is useful in our study because it helps lessen the effects of sampling
errors.  This is an important consideration when investigating
the zeta function because some aspects of its asymptotic behavior 
are reached very slowly (while others are reached extremely rapidly).
In addition, tracking data from different heights (such as $10^{15}$, $10^{19}$, and $10^{22}$) can help us understand how the zeta function approaches its asymptotic behavior.

Ideally one would like to compute the moments over a complete initial interval:

\begin{equation} \label{eq:fstconj}
\frac{1}{T}\int_0^{T}|\zeta(1/2+it)|^{2k}\,dt\,.
\end{equation}

\noindent
However, that could not be done for large values of $T$.  Since we are
interested in asymptotic behavior of the zeta function, we calculated
numerically

\begin{equation} \label{eq:sndconj}
\frac{1}{H}\int_T^{T+H}|\zeta(1/2+it)|^{2k}\,dt\,,
\end{equation}

\noindent
for various choices of $k$, $T,$ and $H$ by evaluating the integral directly.
Section 5 presents the details. We compare the empirical value of (\ref{eq:sndconj}) against 

\begin{equation} \label{eq:shortmoment}
\frac{1}{H}\int_{T}^{T+H} P_k\left(\log \frac{t}{2\pi}\right)\,dt\,, 
\end{equation}

\noindent
where the coefficients of $P_k(x)$ are as given in~\cite{CFKRS1} and~\cite{CFKRS2}. 

We find that for high moments of $|\zeta(1/2+it)|$, and near values of $T$ such 
as $10^{22}$, the variability over adjacent intervals is substantial, even when 
those intervals are long, as long as a block containing $10^9$ zeros; see 
Table~\ref{allsamples}. Moreover, the variability increases rapidly with height, 
a tendency that is observed even for relatively low $T\approx 10^8$ (see the comments following Table~\ref{na1}). 

The leading term prediction (\ref{eq:conjmoments}) does not match well the moment values at the 
heights where we carried our computations. This is not surprising. 
As $k$ increases, the leading coefficient $a(k)g(k)/k^2!$ becomes extremely small (see Table~\ref{akgk}), whereas actual moments of the zeta function do grow moderately rapidly.
Therefore, one cannot hope to get a good approximation for the actual $2k^{th}$ moment from the leading term on its own 
unless $T$ is large enough to enable the leading power $(\log T)^{k^2}$ to compensate for the
small size of the leading coefficient $a(k) g(k)/k^2!$. 
In particular, for high moments (relative to $T$), the contribution of lower order terms is 
likely to dominate.

\begin{table}[ht]
\footnotesize
\caption{\footnotesize  $a(k)$ and $g(k)/(k^2)!$ as given by (\ref{eq:alambda}) and (\ref{eq:cuemoments}) respectively for the first few integers $k$. Numerical values are truncated.}\label{akgk}
\begin{tabular}{|c||c|c|}
\hline
$k$ & $a(k)$ & $g(k)/(k^2)!$ \\
\hline 
1 &  1 &   1\\
2 &  $6.0792 \times 10^{-1}$ & $8.3333 \times 10^{-2}$  \\
3 & $4.9321 \times 10^{-2}$ & $1.1574 \times 10^{-4}$ \\
4 & $2.1468 \times 10^{-4}$ & $1.1482 \times 10^{-9}$  \\
5 &  $3.1326 \times 10^{-8}$ & $4.5202 \times 10^{-17}$ \\
6 &  $1.1415\times 10^{-13}$ & $4.4937 \times 10^{-27}$  \\
7 &  $8.4291 \times 10^{-21}$ & $7.8100 \times 10^{-40}$ \\
8 & $1.0751 \times 10^{-29}$ & $1.7402 \times 10^{-55}$ \\
\hline
\end{tabular}
\end{table}


By contrast, the full moment conjecture of Conrey, Farmer, Keating, Rubinstein, and Snaith ~\cite{CFKRS1}, which gives predictions for lower order terms of $P_k(x)$, and hence takes their contributions into account, is significantly better at matching the values of zeta moments that we computed. For example, using a block of $10^9$  zeros near $T=10^{15}$, and for moments as high as the $12^{th}$ moment, the ratio of (\ref{eq:sndconj}) to (\ref{eq:shortmoment}) differs from the expected 1 by less than $4\%$; see Table~\ref{na3}.  

As remarked earlier, though, our computations do suggest there are large 
oscillations in the remainder terms of the moment conjectures for high moments. 
For example, we find two almost consecutive intervals near $T=10^{22}$, containing 
$10^9$ zeros each, such that the $12^{th}$ moment in one interval is about 16 times 
the $12^{th}$ moment in the other; see Table~\ref{allsamples}. 
Such variability is unlikely to come from
any reasonable conjecture for $P_k(\log t/(2\pi))$.
The reason is that if zeta moments are modeled by such polynomials, then 
perturbing $H$ in (\ref{eq:shortmoment}) is unlikely to produce significant 
enough changes for the value of the integral to change by large multiplicative factors, 
as seen in our data, since $P(\log t/2\pi)$ varies quite slowly as function of $t$.          

To further highlight the difficulties and limitations to be expected in any numerical 
study of high moments, consider Table~\ref{na8}.  It is based on computations
over an interval containing 1.5 million blocks, where each block spans $10^4$ zeros near  $T=10^{22}$.
We find there that more than 50\% of the value of the $12^{th}$ moment comes from just
the 5 blocks that contribute the most.
Interestingly, the size of the $n^{th}$ largest contribution to the $2k^{th}$ moment 
among such blocks is approximated rather well (except for the first few $n$) by a law like

\begin{equation} \label{eq:pl}
\footnotesize
\textrm{$n^{th}$ largest block contribution $\approx$ $n^{-k/5}$ $\times$ (contribution of the largest block).} 
\end{equation}

\noindent
Approximation (\ref{eq:pl}) persists for a long range of $n$; 
see Figure~\ref{powerlaw}.  However, as is explained in Section 4, the constant $5$ 
in this empirical relation is
almost surely a function of the height and number of blocks.

To help understand the observed variability in the moment data, recall that according to a central limit theorem due to Selberg~\cite{S}, the real and imaginary parts of $\log \zeta(1/2+it)/\sqrt{(1/2)\log \log t}$ converge in distribution, over fixed intervals, to 
 independent standard Gaussians. In particular, the distribution of $|\zeta(1/2+it)|$ tends to log-normal. 
So from the outset, we expect the variability in the values of (\ref{eq:sndconj}) to increase 
significantly with $T$ and $k$.  
To iron out such variability, we would like  to take $H$ in integral (\ref{eq:sndconj}) 
relatively large, ideally $H\approx T$, although, from probabilistic considerations,
we might expect that $H\approx T^{1/2}$ might suffice. 
Since this is impractical with current computational resources
when $T= 10^{22}$ say, we typically took $H$ significantly 
smaller than $T$. 
For example, near both $T=10^{15}$ and $T=10^{22}$, the largest we could take $H$ was $\approx 10^8$. 
But then attempts to deduce asymptotics for the integral (\ref{eq:sndconj}) from those 
for the integral (\ref{eq:fstconj}) encounter a problem, which we illustrate in the case of 
the second moment, known to satisfy

\begin{equation} \label{eq:secondmomenttt}
\int_0^T |\zeta(1/2+it)|^2\, dt = T \log \frac{T}{2\pi}+T(2\gamma -1)+E_1(T)\, ,
\end{equation}

\noindent
where $E_1(T)=O(T^{35/108+\epsilon})$; see~\cite{I2}. (This result does not depend on the assumption of the Riemann hypothesis.) Based on (\ref{eq:secondmomenttt}), one might suspect

\begin{equation} \label{eq:shortrangemoment2}
\begin{split}
\frac{1}{H}\int_{T}^{T+H} & |\zeta(1/2+it)|^2\, dt \\
&= \frac{T+H}{H}\log \frac{T+H}{2\pi}-\frac{T}{H}\log \frac{T}{2\pi}+(2\gamma-1)+o(1)\,.
\end{split}
\end{equation}

\noindent
However, in order for  (\ref{eq:shortrangemoment2}) to be valid, we need

\begin{equation} \label{eq:ettunq}
\frac{E_1(T+H)-E_1(T)}{H}=o(1)\,.
\end{equation}

\noindent
Relation (\ref{eq:ettunq}) is certainly true if $H\ge T^{1/2}$ say, 
but not necessarily true if $H$ is substantially smaller. Since in our data the largest $H$ we can take 
is often far below $T^{1/2}$,
the agreement with prediction depends on whether over short intervals the remainder term $E_1(t)$ 
varies slowly enough for something like (\ref{eq:ettunq}) to hold. 
 Of course, a similar statement applies to the remainder terms $E_k(t)$ corresponding to higher 
moments. Therefore, in summary, the variability of the moment data 
 is a function of the following:




\begin{itemize}
\item
is the height $T$ large enough for asymptotics to apply? 
\item
are the number and size of samples large enough to be representative of a log-normal variable? 
\item
given $k$, are there significant oscillations in the remainder term $E_k(T)$? 
\end{itemize}

Numerical data show some long-range correlations (on the scale of a few thousand zeros)
for the second moment; see Figure~\ref{blockcorr}.
However, we do not find similar evidence for such correlations in the case of higher moments. 
To further examine the correlations in the second moment, we numerically investigate the 
shifted fourth moment,
where K\"osters~\cite{Ko} proved a kernel law for shifts on the scale of mean spacing of zeros.
For larger shifts, we observe a departure from this law,
and the onset of asymptotic oscillations; see Figures~\ref{smg1} and~\ref{smg2}.

In the course of performing our moment computations, various local models of the zeta function (i.e. models to numerically approximate $|\zeta(1/2+it)|$ over short intervals) were analyzed. The results of those analyses are of independent interest, for guidance in
selection of numerical methods for other calculations that might be done in
the future with the zeta function. An attractive local model that was tested is due to 
Gonek, Hughes, and Keating~\cite{GHK}.  Another model, which numerically works quite well, 
is to approximate $|\zeta(1/2+it)|$ by a suitably normalized polynomial that vanishes 
at the zeta zeros near $t$. Numerical
computations suggest this approximation converges linearly in the number of zeros 
used.  (We consequently found an elementary proof of this.)

The paper is organized as follows. In Section 2, we provide a further discussion of 
some known results and conjectures for the moments of $|\zeta(1/2+it)|$. In Section 3 we document the 
datasets used in our study, and calculate the various moment predictions for them. 
Section 4 contains the results of our numerical computations of the $2^{nd}$ to $12^{th}$ 
 even moments. In Section 5, we outline the numerical methods employed, which include point-wise approximations to zeta (the main one being the Odlyzko-Sch\"onhage 
algorithm), the choice of the integration technique, 
and implementation details. In Section 6, we numerically verify the accuracy of some local 
(short-range) point-wise approximations of zeta that were used in our computations.
The results of these computations are of independent interest.

\section{Results and conjectures for zeta moments}




Number-theoretic methods were not able to predict, in general, the value of the factor $g(k)$, which occurs in asymptotic (\ref{eq:conjmoments}). 
In contrast, a general prediction for the arithmetic factor $a(k)$ was established explicitly as

\begin{equation} \label{eq:alambda}
a(k):=\prod_p\left((1-1/p)^{k^2}\left(\sum_{m=0}^\infty d_{k}^2(p^m)p^{-m}\right)\right)\,,
\end{equation}

\noindent
where $d_{k}$ is the $k$-fold divisor function. We remark the size of $a(k)$ is mostly determined by the contributions of the first $O(k^2)$ primes. 



Although $g(k)$ is not critical to determining the order of \mbox{$\zeta(1/2+it)$}, 
 it is important in comparing empirical data to conjectures due to the extra precision it provides.
 A conjecture for the value 
of $g(k)$ was first provided by Keating and Snaith~\cite{KS}.
Motivated by recent progress in developing (conjectural, numerical, and heuristic)
connections between RMT and the zeta function, they suggested there
might be a relation between moments of the characteristic polynomials of unitary matrices
and moments of the zeta function. This led them to conjecture:

\begin{equation} \label{eq:cuemoments}
\frac{g(k)}{k^2!}:=\prod_{j=0}^{k-1} \frac{j!}{(j+k)!} \,.
\end{equation}

Specifically, Keating and Snaith considered $Z_N(A,\theta)$, the characteristic polynomial (in $e^{i\theta}$) of an \mbox{$N\times N$} unitary matrix $A$. 
They used Selberg's formula to calculate the expected moments of $|Z_N(A,\theta)|$, where the expectation is taken with respect to the normalized Haar measure on the group of $N\times N$ unitary matrices. They found for integer values of $k >0$ that

\begin{equation} \label{eq:cuemoments0}
\mathbb{E}_N |Z_N(A,\theta)|^{2k}=\prod_{j=0}^{N-1} \frac{j!\,(j+2k)!}{(j+k)^2!}\,. 
\end{equation}

\noindent
The right side of formula (\ref{eq:cuemoments0}) is a polynomial in $N$ of degree $k^2$, with leading coefficient given by the right side of (\ref{eq:cuemoments}). (The full result of~\cite{KS} is more general as it can be continued to the half-plane $\Re(k)>-1/2$.)



The absence of the arithmetic factor $a(k)$ from formula (\ref{eq:cuemoments0}) hints that the moments of zeta
might split as the product of two parts, one that is universal and is due to RMT
fluctuations (on the scale of mean spacing), and another that is specific to zeta and corresponds to the contributions of small primes; see~\cite{BK} and~\cite{GHK}. A similar phenomenon appears in some statistics of zeta zeros, such as zero correlations and number-variance of zero spacings.

It appears, however, that RMT is only able to predict the universal part of 
the leading term asymptotic for zeta moments. 
For example, Ivi\'c~\cite{I1} and, separately, Conrey~\cite{C}, determined explicitly the coefficients of the fully proven fourth moment polynomial $P_2(x)$. Their results lead to

\begin{equation} \label{eq:exp4}
\begin{split}
\frac{1}{T}\int_0^T |\zeta(1/2+it)|^4\, dt \approx\, &\, 0.050660 \left(\log \frac{T}{2\pi}\right)^4 +\, 0.496227 \left(\log \frac{T}{2\pi}\right)^3\\
&+\, 0.937279\left(\log \frac{T}{2\pi}\right)^2+\,1.35334\left(\log \frac{T}{2\pi}\right)\\
&-\,0.040924  +E_2(T)/T\,,
\end{split}
\end{equation}
 
\noindent
where numerical values are obtained via truncation. 
By comparison, the straightforward RMT-based prediction for the fourth moment is determined by a quantity like

\begin{equation} \label{eq:rmt1}
a(2)\, \mathbb{E}_N |Z_N(A,\theta)|^{4}=a(2)\, \prod_{j=1}^N\frac{\Gamma(j)\Gamma(j+4)}{(\Gamma(j+2))^2}\, ,   
\end{equation}

\noindent
where $a(2)=0.607927\ldots\,\,$. If expression (\ref{eq:rmt1}) correctly predicts lower order terms for the fourth moment, then it should yield  something similar to right side of (\ref{eq:exp4}), but it does not; 
if we make the identification $N\to \log T/(2\pi)$, as suggested by~\cite{KS}, then expression (\ref{eq:rmt1}) produces the following polynomial instead:

\begin{equation} \label{eq:exprmt4}
\begin{split}
0.050660\left(\log \frac{T}{2\pi}\right)^4 +\, 0.405284 & \left(\log \frac{T}{2\pi}\right)^3 +1.16519\left(\log \frac{T}{2\pi}\right)^2 \\
&+\,1.41849\left(\log \frac{T}{2\pi}\right)+\, 0.607927\,.
\end{split}
\end{equation}

\noindent
So, the leading terms in polynomials (\ref{eq:exp4}) and (\ref{eq:exprmt4}) agree, but lower order terms do not. 
This discrepancy indicates RMT does not model zeta moments well enough to enable correct predictions of 
lower asymptotic terms.  However, in numerical studies of zeta moments, it is useful to go beyond the 
leading coefficient, to have predictions for lower order coefficients in the moment polynomial $P_k(x)$ 
in conjecture (\ref{eq:rubmomentold10}).


This is one of the reasons the recent work of Conrey, Farmer, Keating, Rubinstein, and
Snaith~\cite{CFKRS1}~\cite{CFKRS2} has been useful to our numerical investigation.
They gave a conjecture for $P_k(x)$ which agrees with the known cases $k=1$ and $k=2$, as well as with the RMT 
prediction (\ref{eq:conjmoments}). 
Their conjecture is moreover supported by empirical data at low heights (near $T= 2\times 10^6$). 
Importantly,  they computed the coefficients of $P_k(x)$ for the first few integers $k$. 
This enabled us to calculate their predictions for moments beyond the fourth moment ($2k=4$).

\section{Preliminaries}

We sample consecutive blocks of zeros $B_n$. The blocks $B_n$ are of equal size, and are located in a neighborhood of the height $T$. We denote the size of a block $B_n$  by $|B_n|$, and choose $|B_n|=|B_{n+1}|$ for all relevant $n$. We consider many ratios of the form

\begin{equation}\label{eq:appmo2}
\frac{\int_{ B_n}|\zeta(1/2+it)|^{2k}\,dt}{\int_{B_n} P_k\left(\log \frac{t}{2\pi}\right)\,dt}\,, 
\end{equation}

\noindent
where, by an abuse of notation, the symbolism $\int_{B_n}$ is short for integrating over the interval spanned by block $B_n$. So if $\alpha_n$ and $\beta_n$ denote the ordinates of the first and last zeros in block $B_n$, respectively, then $\int_{B_n}:=\int_{\alpha}^{\beta}$.  We expect the average of many ratios of the form (\ref{eq:appmo2}) to approach 1. Notice if $T$ is large, and the length of the interval spanned by block $B_n$ is small compared to $T$, then the denominator in ratio (\ref{eq:appmo2}) is largely a function of $T$ multiplied by the length of the interval spanned by $B_n$.  

 In our computations, we aggregated the moment data for each 1,000 consecutive zeros.  (Thus blocks
did not have the same lengths.  However, since zeros are spaced quite regularly, the differences
in lengths of blocks with the same number of zeros at the same height were minor.) 
The block sizes used in our computations range anywhere from $10^3$ to $10^9$ consecutive zeros (data for a block size of $10^4$, for instance, was obtained by averaging the aggregated data for 10 successive blocks of 1,000 zeros each); see Section 5 for details. Most of our computations were 
performed in the vicinity of the $10^{23}$-rd zero, which is near $t=1.3\times 10^{22}$. 

To investigate the behavior of the remainder term $E_k(t)$ in the full moment conjecture~\cite{CFKRS1}, we numerically examine the quantity

\begin{equation}\label{eq:apper2}
\log\left(\frac{\int_{ B_n}|\zeta(1/2+it)|^{2k}\,dt}{\int_{B_n} P_k\left(\log \frac{t}{2\pi}\right)\,dt}\right) = \log\left(1+\frac{E_k(\beta_n)-E_k(\alpha_n)}{\int_{\alpha_n}^{\beta_n} P_k\left(\log \frac{t}{2\pi}\right)\,dt}\right)\,, 
\end{equation}

\noindent
where $\alpha_n$ and $\beta_n$ denote the ordinates of the first and last zeros of block $B_n$, respectively. 

Table~\ref{ta1} contains the exact coordinates of the sets for which we computed moments. 
Except for set $S8$, all these zeros were computed by the second
author at AT\&T Labs--Research in the late 1990s, in an extension of the earlier
computations described in \cite{O1}, and will be documented in the book
\cite{HO}. Some output files from those computations 
were corrupted by buggy archiving software during their transfer from AT\&T.
Although we managed to clean up many of these files, some data was lost irretrievably. 
As a result, there were a few instances where we have missing blocks of
anywhere between 1,000 to 4,000 consecutive zeros.  In such cases,  we skipped
the missing blocks.

\begin{table}[ht] 
\footnotesize
\caption{\footnotesize The datasets}\label{ta1}
\begin{tabular}{|c|c|l|}
\hline
Set & Approximate span & First and last zeros, respectively \\
\hline
$A23$ & 5 billion zeros & $1.30664344087953251142539323425414 \times 10^{22}$ \\
&& $1.30664344087959822199974045053551\times 10^{22}$\\
\hline
$B23$ & 11 billion zeros & $1.30664344087935097997293481220857 \times 10^{22}$\\
&&  $1.30664344087949333176034585636412\times 10^{22}$  \\
\hline
$O20$ &1 billion zeros & $1.52024401160089830109496959179 \times 10^{19}$ \\
&& $1.52024401161628795187223388010 \times 10^{19}$ \\
\hline
$Z16$ & 1 billion zeros & $2.5132741232472002749333722 \times 10^{15}$ \\
&&$2.5132743103949376298283407 \times 10^{15}$ \\ 
\hline
$S8$ & 100 million zeros &  14.1347251417347 \\
&&  42653549.7609516\\
\hline
\end{tabular}
\end{table}

Table~\ref{tabexpmo} contains the moment values predicted by the leading term (\ref{eq:conjmoments}) at $T=10^{22}$. Predictions start to grow more slowly at the $14^{th}$ moment and they 
completely collapse by the $28^{th}$ moment. But by the asymptotic relation (\ref{eq:secondmomenttt}), we have

\begin{equation} \label{eq:incmots1}
 \frac{1}{T}\int_0^T|\zeta(1/2+it)|^{2}\,dt \ge 1\,. 
\end{equation}

\noindent
for $T>T_0$, where $T_0>0$ is some absolute number.  
It follows by induction and Jensen's inequality that for $k \ge 1$ and $T>T_0$ 

\begin{equation} \label{eq:incmots}
1\le \frac{1}{T}\int_0^T|\zeta(1/2+it)|^{2k}\,dt \le \frac{1}{T}\int_0^T|\zeta(1/2+it)|^{2k+2}\,dt\,.
\end{equation}

\noindent
Therefore, the $2k^{th}$ moments should increase with $k$ for $T>T_0$. 
 In view of Table~\ref{tabexpmo}, the height $T$ must then be larger than $10^{22}$ in order to reach the 
 leading term asymptotic of certainly the $18^{th}$ moment, and most likely the $14^{th}$ and higher 
moments. This is the reason we chose to cut off our computation at the $12^{th}$ moment.

\begin{table}[ht]
\footnotesize
\caption{\footnotesize The expected moments as given by the right side of (\ref{eq:conjmoments}) for $T$ equal to the ordinate of the first zero in set $A23$.\label{tabexpmo}}
\begin{tabular}{|c|c|}
\hline
Moment & Value \\
\hline
2 & 5.09 $\times 10^1$ \\ 
4 & 3.40 $\times 10^5$ \\
6 & 1.31 $\times  10^{10}$ \\
8 & 5.04 $\times 10^{14}$ \\
10 & 6.67 $\times 10^{18}$ \\
12 & 1.44 $\times 10^{22}$ \\
14 & 2.86 $\times 10^{24}$ \\
16 & 3.27 $\times 10^{25}$ \\
18 & 1.44 $\times 10^{25}$ \\
20 & 1.74 $\times 10^{23}$ \\
22 & 4.29 $\times 10^{19}$ \\
24 & 1.64 $\times 10^{14}$ \\
26 & 7.61 $\times 10^6$ \\
28 & 3.50 $\times 10^{-3}$ \\
\hline
\end{tabular}
\end{table}

Table~\ref{exp1} contains the full moment predictions (\ref{eq:shortmoment}) for the various sets listed in Table~\ref{ta1}, where the coefficients of the polynomial $P_k(x)$ in (\ref{eq:shortmoment}) are as computed in~\cite{CFKRS1} and~\cite{CFKRS2}. 
Table~\ref{exp2} contains the moment predictions when the polynomial $P_k(x)$ in (\ref{eq:shortmoment}) is exchanged for its leading term only.

\begin{table}[ht]
\footnotesize
\caption{\footnotesize The expected $2k^{th}$ moment (truncated here after three digits) as given by integral (\ref{eq:shortmoment}), where the coefficients of $P_k(x)$ are as computed in~\cite{CFKRS1} and~\cite{CFKRS2}. For each set, the integral is evaluated by setting $T$ to the ordinate of the first zero, and $T+H$ to the ordinate of the last zero. The first and last zeros of each set are specified in Table~\ref{ta1}.\label{exp1}}
\begin{tabular}{|c|c|c|c|c|c|c|}
\hline
Set & $2k=2$  & $2k=4$ & $2k=6$ & $2k=8$ & $2k=10$ & $2k=12$ \\
\hline
$A23$ & $5.02 \times 10^1$  &  $3.82\times 10^5$ & $3.30 \times 10^{10}$ & $1.04 \times 10^{16}$ &  $7.32 \times 10^{21}$ &  $8.89 \times 10^{27}$ \\
$B23$ & $5.02 \times 10^1$  &  $3.82\times 10^5$ & $3.30 \times 10^{10}$ & $1.04 \times 10^{16}$ &  $7.32 \times 10^{21}$ &  $8.89 \times 10^{27}$ \\
$O20$ & $4.34 \times 10^1$ & $2.20 \times 10^5$ & $1.03 \times 10^{10}$  & $1.54 \times 10^{15}$ & $4.66 \times 10^{20}$ &  $2.26 \times 10^{26}$ \\
$Z16$ & $3.47 \times 10^1$ & $9.41 \times 10^4$ & $1.77 \times 10^9$  & $8.77 \times 10^{13}$ & $7.63 \times 10^{18}$ &  $9.70 \times 10^{23}$ \\
$S8$ & $1.58 \times 10^1$ & $5.28 \times 10^3$ & $5.58 \times 10^6$  & $9.87 \times 10^{9}$ & $2.33 \times 10^{13}$ &  $6.67 \times 10^{16}$ \\
\hline
\end{tabular}
\end{table}

\begin{table}[ht]
\footnotesize
\caption{\footnotesize The expected $2k^{th}$ moment (truncated here after three digits) when the polynomial $P_k(x)$ in (\ref{eq:shortmoment}) is replaced by its leading term only, which has coefficient $a(k)g(k)/k^2!$. For each set, the integral is evaluated by setting $T$ to the ordinate of the first zero, and $T+H$ to the ordinate of the last zero. The first and last zeros of each set are specified in Table~\ref{ta1}.\label{exp2}}
\begin{tabular}{|c|c|c|c|c|c|c|}
\hline
Set & $2k=2$  & $2k=4$ & $2k=6$ & $2k=8$ & $2k=10$ & $2k=12$ \\
\hline
$A23$ & $4.90 \times 10^1$ & $2.94 \times 10^5$ & $9.44 \times 10^9$  & $2.80 \times 10^{14}$ & $2.66 \times 10^{18}$ &  $3.84 \times 10^{21}$ \\
$B23$ & $4.90 \times 10^1$ & $2.94 \times 10^5$ & $9.44 \times 10^9$  & $2.80 \times 10^{14}$ & $2.66 \times 10^{18}$ &  $3.84 \times 10^{21}$ \\
$O20$ & $4.23 \times 10^1$ & $1.62 \times 10^5$ & $2.49 \times 10^9$  & $2.61 \times 10^{13}$ & $6.56 \times 10^{16}$ &  $1.85 \times 10^{19}$ \\
$Z16$ & $3.36 \times 10^1$ & $6.47 \times 10^4$ & $3.13 \times 10^8$  & $6.57 \times 10^{11}$ & $2.07 \times 10^{14}$ &  $4.66 \times 10^{15}$ \\
\hline
\end{tabular}
\end{table}

The predictions (\ref{eq:shortmoment}) for the sets $A23$, $B23$, $O20$, $Z16$, and $S8$
are calculated by evaluating the integral (\ref{eq:shortmoment}) for the exact coordinates given in Table~\ref{ta1}. Notice when $T$ is large, predictions are insensitive to the precise choice of $H$ as long as $H$ is not too large, say $H<T^{1/2}$. This will be the case for most samples from those sets. 

Finally, we point out the following notational convention in the next section. The heading of many of the tables in Section 4 will have the format 

\begin{equation}
\textrm{$X$ samples $Y$}  \nonumber
\end{equation}

\noindent
This means the values listed in the table are based on $X$ (usually consecutive) blocks, where each block consists of $Y$ consecutive zeros. For simplicity, the block size $Y$ is often written in abbreviated form; e.g. if $Y$ is 1M, this means each block consists of a million consecutive zeros. Similarly, 1B means each block consists of a billion consecutive zeros, and so on. We frequently refer to values of the moments over blocks simply as samples.

\section{Numerical results}

At the heights used in our study, there are substantial disparities between the full moment 
prediction (\ref{eq:shortmoment}) and  the leading term prediction (obtained by 
replacing $P_k(x)$ in (\ref{eq:shortmoment}) by its leading terms only).  
So, a relatively a good agreement with one conjecture implies lack of agreement with the other. 
For this reason, let us focus our attention on the full moment prediction (\ref{eq:shortmoment}), which is believed to be more accurate. 

Conrey, Farmer, Keating, Rubinstein, and Snaith~\cite{CFKRS1} (see also~\cite{CFKRS2}) considered ratios of the form 

\begin{equation} \label{eq:appmo1}
\frac{\int_0^{T}|\zeta(1/2+it)|^{2k}\,dt}{\int_0^{T} P_k\left(\log \frac{t}{2\pi}\right)\,dt}\,. 
\end{equation}

\noindent
for $T$ near $10^6$. They found a very good agreement between the data and the full moment predictions at that height. We remark the predictions for high moments near $T=10^6$, and also near $T=10^7, 10^{15}, 10^{19}, 10^{22}$ (heights used in our experiments), are not really determined by the leading term asymptotic, but by the lower order terms in the full moment conjecture. This is because lower order terms still contribute the most even for $T$ as large as $10^{22}$.

We first discuss the moment data at low heights near $10^7$. This is set $S8$ in Table~\ref{ta1}, which consists of the first $10^8$ zeros. To aid in the analysis, let us define the following subsets: 

\begin{table}[ht]
\footnotesize
\caption{\footnotesize Some subsets in $S8$. \label{tabs81}} 
\begin{tabular}{|c|l|l|l|}
\hline
Set & Initial zero & Final zero & Size of the set \\ 
\hline
$s1$ &  10,000,000 & 10,100,000  & $10^5$ zeros\\ 
$s2$ & 90,000,000 & 90,100,000  & $10^5$ zeros \\
$s3$ &  10,000,000 & 11,000,000 & $10^6$ zeros\\ 
$s4$ &  90,000,000 & 91,000,000 & $10^6$ zeros \\
$s5$ &  10,000,000 & 20,000,000 & $10^7$ zeros\\ 
$s6$ & 90,000,000 & 99,980,000 & $0.9998 \times 10^7$ zeros \\
\hline
\end{tabular}
\end{table}

Table~\ref{na1} lists the value of the ratio (\ref{eq:appmo2}) for each of the subsets defined in Table~\ref{tabs81}. For example, the first entry in Table~\ref{na1}, which corresponds to $2k=2$ and subset $s1$, was calculated using the formula:

\begin{equation}
\frac{\int_{s_1}|\zeta(1/2+it)|^{2}\,dt}{\int_{s_1} P_1\left(\log \frac{t}{2\pi}\right)\,dt} = \frac{\int_{\gamma_{10^7}}^{\gamma_{10^7+10^5}}|\zeta(1/2+it)|^{2}\,dt}{\int_{\gamma_{10^7}}^{\gamma_{10^7+10^5}} P_1\left(\log \frac{t}{2\pi}\right)\,dt}\,.\ 
\end{equation}

\begin{table}[htbp]
\footnotesize
\caption{\footnotesize Ratio of empirical moment to full moment prediction for the zero subsets defined in Table~\ref{tabs81}. 
\label{na1}}
\begin{tabular}{|c||c|c||c|c||c|c||}
\hline
$2k$ & $s1 $ & $s2$ & $s3$ & $s4$ & $s5$ & $s6$  \\ 
\hline
2 &  1.000 & 1.001 & 1.000 &  1.000 & 1.000 & 1.000  \\ 
4 &  0.996 & 1.007 & 1.000 &  1.000 & 1.001 & 1.000 \\ 
6 &  0.975  & 0.999 & 0.997 & 0.996 &  1.001  & 1.000 \\ 
8 &  0.943  & 0.981 & 0.992 & 0.983 &  1.001 & 0.999 \\ 
10 & 0.909 & 0.962 &  0.989 & 0.962 & 1.001 & 1.000 \\ 
12&   0.875 & 0.940 & 0.987 & 0.931 & 1.000  & 1.004 \\ 
\hline
\end{tabular}
\end{table}

Table~\ref{na1} indicates a fairly good agreement with the full moment prediction (\ref{eq:shortmoment}). However, even at these modest heights, we find evidence for substantial variability in the moment data. For instance, we can find blocks of $10^5$ zeros each such that the $12^{th}$ moment is half what is predicted by (\ref{eq:shortmoment}) (e.g. block $[\gamma_{5\times 10^7},\gamma_{5\times 10^7+10^5}]$). Also, the variability in the moment data seems to increase with height. 
For example, when we sample 100 consecutive blocks of $10^5$ zeros each near zero numbers 21 million 
and 91 million, the standard deviation of the $12^{th}$ moment, divided by the full moment 
prediction (\ref{eq:shortmoment}), increases from about $0.994$ to about $2.385$, which seems 
substantial considering the statistics discussed here change on a logarithmic scale.

We next consider the situation higher on the critical line. 
Table~\ref{na3} contains statistical summaries for moment values from set $Z16$, 
which is a set of $10^9$ zeros near the $10^{16}$-th zero.
The first five columns of Table~\ref{na3} were calculated as follows. 
We subdivided set $Z16$ 
into $10^6$ blocks $B_n$, where each block has $10^3$ consecutive zeros.
(As was mentioned earlier, we aggregated the moment data for each 1000 consecutive zeros, 
so the smallest block size we can work with is 1000; see Section 5 for details.) 
Next, for each $k\in \{1,\ldots,6\}$, we computed the quantities

\begin{equation}
I_n:=\frac{1}{\beta_n-\alpha_n}\int_{\alpha_n}^{\beta_n}|\zeta(1/2+it)|^{2k}\,dt\,,\qquad n\in [1,10^6]\,,\
\end{equation}

\noindent
where $\beta_n$ and $\alpha_n$ are the ordinates of the first and last zeros of block $B_n$ respectively (the blocks are chosen so the last zero of $B_n$ is the first zero of $B_{n+1}$). Lastly, for each $1\le j\le 1000$ we calculated 

\begin{equation} \label{eq:sndm}
x_j:=\frac{\frac{1}{1000}\sum_{n=1+(j-1)1000}^{j1000} I_n}{\frac{1}{H}\int_T^{T+H} P_k(t)\,dt}\,,
\end{equation}
 
\noindent
where $T$ and $T+H$ are chosen equal to the ordinates of the first and last zeros of the set $Z16$, which is where the blocks $B_n$ reside. Notice the denominator in (\ref{eq:sndm}) remains practically constant for $H$ small compared to $T$, and so it is essentially the same as writing 

\begin{equation}
\frac{1}{\beta_{1000 j}-\alpha_{1+1000(j-1)}}\int_{\alpha_{1+1000(j-1)}}^{\beta_{1000j}} P_k\left(\log \frac{t}{2\pi}\right)\,dt\,,
\end{equation}

\noindent
where $\alpha_{1+1000(j-1)}$ is the ordinate of the first zero of $B_{1+1000(j-1)}$, and $\beta_{1000j}$ is the ordinate of the last zero of $B_{1000j}$.  This procedure yields 1000 numerical values (i.e. sample points) $\{x_j, 1\le j\le 1000\}$, where each $x_j$ now corresponds to a block of $10^6$ zeros.  The first five columns of Table~\ref{na3} contain the mean, minimum, and maximum of these values, as well as their standard deviation about the empirical mean. The other columns in the table were calculated similarly, but using larger block sizes ($10^7$, and $10^8$, zeros, respectively).

\begin{table}[ht]
\footnotesize
\caption{\footnotesize Ratio of empirical moment to full moment prediction for samples from $Z16$. The full moment prediction is given in Table~\ref{exp1}. The columns Min, Max, Mean, and SD refer to the min, max, mean, and the standard deviation of the ratios about their empirical mean.}\label{na3}
\begin{tabular}{|c|c||c|c|c||c|c|c||c|c|c||}
\hline$2k$ & Mean &\multicolumn{3}{|c||}{$Z16$: 1000 samples 1M} & \multicolumn{3}{|c||}{$Z16$: 100 samples 10M} & \multicolumn{3}{|c||}{$Z16$: 10 samples 100M} \\
\cline{3-11}
 &   &Min  & Max & SD&  Min & Max & SD& Min & Max & SD\\\hline
2 & \textbf{1.000} & 0.987 & 1.009  & \textbf{0.003}& 0.999 & 1.002 & \textbf{0.001}& 1.000 & 1.000 & \textbf{0.000}\\ 
4 & \textbf{1.000} & 0.831 & 1.602 & \textbf{0.077}& 0.969 & 1.048 & \textbf{0.017} & 0.996 & 1.006 & \textbf{0.003}\\ 
6 & \textbf{1.003} & 0.494 & 8.687 & \textbf{0.530}& 0.805 & 1.796 & \textbf{0.143}& 0.966 & 1.079 & \textbf{0.035}\\8 & \textbf{1.019} & 0.178 & 37.68  & \textbf{1.837}& 0.532 & 4.914 & \textbf{0.555}& 0.882 & 1.262 & \textbf{0.143}\\ 
10 & \textbf{1.039} & 0.046 & 101.0 & \textbf{4.202} & 0.277 & 11.57 & \textbf{1.327} & 0.709 & 1.800 & \textbf{0.354}\\
12 & \textbf{1.035} & 0.009 & 190.4 & \textbf{7.182} & 0.116 & 20.67 & \textbf{2.301}& 0.484 & 2.553 & \textbf{0.628}\\
\hline
\end{tabular}
\end{table}


In view of Table~\ref{na3}, we see the standard deviation of samples generally declines 
like $1/\sqrt{|B|}$, where $|B|$ is the  block size. This indicates that the moments 
over such long blocks (millions of zeros each) act statistically independently. 
Also, the range of values spanned by high moments (i.e. the interval [Min, Max]) appears to 
decline linearly with $|B|$. This observation is likely  due to the sparsity of blocks with 
large contributions. For example, if the block size is increased from 1M to 10M, then since 
blocks with large contributions are rare, and since they usually do not occur near each 
other (at least as far as our data is concerned), the range [Min, Max] should shrink to 
something like [Min, Max/10], as observed.



Table~\ref{na3} shows reasonable agreement with the full moment 
prediction (\ref{eq:shortmoment}) near $T=10^{15}$.
At larger heights, the agreement is poorer.
This is especially true for high moments, roughly the $8^{th}$ and higher moments. 
For instance, consider Table~\ref{na4}, which is based on a set of $1.5 \times 10^{10}$ 
zeros in the vicinity of the $10^{23}$-rd zero. There, the $12^{th}$ moment is about half what is expected. 
This growing disparity between prediction and evidence is likely due to a substantial
extent to the fact that the integration interval is considerably shorter relative
to the height around the $10^{23}$-rd zero.

\begin{table}[ht]
\footnotesize
\caption{\footnotesize Ratio of empirical moment to full moment prediction for samples from $A23\cup B23$ (full moment prediction is given in Table~\ref{exp1}). The columns Min, Max, Mean, and SD refer to the min, max, mean, and the the standard deviation of the ratios about their empirical mean. }\label{na4}
\begin{tabular}{|c|c||c|c|c||c|c|c||}
\hline
$2k$ & Mean &\multicolumn{3}{|c||}{$A23 \cup B23$: 15,000 samples 1M} & \multicolumn{3}{|c||}{$A23 \cup B23$: 1500 samples 10M}  \\
\cline{3-8}
 &   & $\,\,\,$ Min $\,\,\,$ & $\,\,\,$ Max $\,\,\,$  & SD&  $\,\,\,$ Min $\,\,\,$ & $\,\,\,$ Max $\,\,\,$ & SD\\ 
\hline
2 & \textbf{1.000} & 0.978 & 1.045  &  \textbf{0.007} & 0.993 & 1.008 & \textbf{0.002}  \\ 
4 & \textbf{1.000} & 0.642 & 8.472  &  \textbf{0.270} & 0.859 & 1.895 & \textbf{0.082} \\ 
6 & \textbf{0.990} & 0.172 & 121.4  &  \textbf{2.289} & 0.468 & 13.38 & \textbf{0.735} \\  
8 & \textbf{0.935} & 0.021 & 556.1  &  \textbf{8.179} & 0.142 & 56.51 & \textbf{2.623} \\  
10 & \textbf{0.791} & 0.001 & 1184  &  \textbf{15.52} & 0.025 & 118.8 & \textbf{4.945} \\  
12 & \textbf{0.574} & 0.000 & 1486  &  \textbf{18.27} & 0.003 & 148.7 & \textbf{5.799} \\  
\hline
\end{tabular}
\end{table}

Table~\ref{na5} provides a summary of all our moment data. 
It shows that at large heights (e.g. near $T=10^{22}$), even a block size of $10^9$ 
zeros is not enough to control the variability in the moment data over adjacent blocks.
An extreme example is that of the adjacent blocks $b5$ and $b6$, where the $12^{th}$ 
moment is multiplied by more than a factor of 16 from one block to the next.

\begin{table}[ht]
\footnotesize
\caption{\footnotesize Ratio: empirical moment/full moment prediction, for subsets from $S8$, $Z16$, $O20$, $A23$ and $B23$.  Each subset consists of $10^9$ consecutive zeros except for $s8$, which consists of about $4\times 10^7$ zeros (specifically, zeros 60M to 99.98M in $S8$). Subsets $ax$ and $bx$ are of increasing height, and are approximately consecutive except for some small gaps. The column $\tilde{T}$ is the approximate height at which the subset starts, and column $\tilde{H}$ is the approximate length of the interval spanned by the subset.\label{allsamples}}\label{na5}
\begin{tabular}{|c||c|c|c|c|c|c||c|c|}
\hline
Sample & $2k=2$ & $2k=4$ & $2k=6$ & $2k=8$ & $2k=10$ & $2k=12$ & $\tilde{T}$ & $\tilde{H}$ \\
\hline
$s8$  & 1.000 & 1.000 & 1.000 & 1.001 & 1.002 & 1.004 & $10^7$ & $10^7$ \\
$z16$ & 1.000 & 1.000 & 1.003 & 1.019 & 1.039 & 1.035 & $10^{15}$ & $10^8$\\
$o20$ & 1.000 & 1.000 & 1.003 & 0.982 & 0.873 & 0.667 & $10^{19}$ & $10^8$ \\
$b1$ & 1.000 & 0.997 & 0.941 & 0.731 & 0.419 & 0.176 & $10^{22}$ & $10^8$ \\ 
$b2$ & 1.000 & 1.006 & 1.067 & 1.183 & 1.166 & 0.908 & $10^{22}$ & $10^8$\\ 
$b3$ & 1.000 & 1.002 & 0.989 & 0.864 & 0.587 & 0.297& $10^{22}$ & $10^8$\\ 
$b4$ & 1.000 & 0.994 & 0.931 & 0.740 & 0.454 & 0.208& $10^{22}$ & $10^8$\\ 
$b5$ & 1.000 & 0.991 & 0.895 & 0.632 & 0.324 & 0.123& $10^{22}$ & $10^8$\\ 
$b6$ & 1.000 & 1.008 & 1.127 & 1.543 & 2.011 & 2.011& $10^{22}$ & $10^8$\\ 
$b7$ & 1.000 & 0.991 & 0.932 & 0.772 & 0.517 & 0.269& $10^{22}$ & $10^8$\\ 
$b8$ & 1.000 & 0.999 & 0.974 & 0.836 & 0.570 & 0.303& $10^{22}$ & $10^8$\\ 
$b9$ & 1.000 & 1.003 & 0.980 & 0.835 & 0.568 & 0.306& $10^{22}$ & $10^8$\\ 
$b10$ & 1.000 & 0.988 & 0.903 & 0.685 & 0.399 & 0.174& $10^{22}$ & $10^8$\\ 
$a1$ & 1.000 & 0.995 & 0.922 & 0.684 & 0.374 & 0.151& $10^{22}$ & $10^8$\\ 
$a2$ & 1.000 & 1.002 & 1.005 & 0.976 & 0.820 & 0.542& $10^{22}$ & $10^8$\\ 
$a3$ & 1.000 & 1.005 & 1.108 & 1.449 & 1.820 & 1.823& $10^{22}$ & $10^8$\\ 
$a4$ & 1.000 & 1.006 & 1.061 & 1.185 & 1.213 & 1.003& $10^{22}$ & $10^8$\\ 
$a5$ & 1.000 & 1.007 & 1.019 & 0.911 & 0.626 & 0.321& $10^{22}$ & $10^8$\\ 
\hline
\end{tabular}
\end{table}

Since high moments are determined by a few samples with large contributions, it is useful to measure the impact of such samples more accurately. So, consider Table~\ref{na7}, which is a ``moments of the moments'' table calculated near zero number $10^{23}$.

\begin{table}[ht]
\footnotesize
\caption{\footnotesize Moments of the ratios empirical moment/predicted moment after being normalized to have mean 0 and variance 1 (predicted moment is given in Table~\ref{exp1}). The columns $p=3$, $p=4$,...etc  refer to the third moment, fourth moment,...etc, of the quantity empirical moment/predicted moment.}\label{na7}
\begin{tabular}{|c|l|l|l|l|l|l|}
\hline
$2k$ &\multicolumn{6}{|c|}{$A23 \cup B23$: 150,000 samples 100,000} \\
\cline{2-7}
 & $p=3$ & $p=4$ & $p=5$ & $p=6$ & $p = 7$ & $p=8$\\ 
\hline
2 & 1.051 & 7.251 & 42.41 & 411.7 & 4664 & 58750\\ 
4 & 20.61 & 1042 & 69580 & 5135000 & 396100000 & 31280000000\\ 
6 & 83.86 & 10960 & 1599000 & 243800000 & 37990000000 & 6000000000000\\ 
8 & 146.6 & 26920 & 5264000 & 1059000000 & 216300000000 & 44580000000000\\ 
10 & 185.6 & 39590 & 8836000 & 2014000000 & 464100000000 & 107700000000000\\ 
12 & 209.5 & 48560 & 11660000 & 2845000000 & 700300000000 & 173400000000000\\ 
\hline
\end{tabular}
\end{table}

In Table~\ref{na7}, there is a notable slowdown in the growth rate with $k$ of the ``moments of the moments.'' This slowdown is directly related to the frequency and precise size of large block contributions. To get a better sense of the distribution of such contributions, let us consider Table~\ref{na8}. It lists the sum of the $n$ largest block contributions to each moment for several $n$. To explain how the table was constructed, take the entries corresponding to $2k=2$ and $n= 1,2,\ldots,5$ for example. They were calculated by first computing $\{x_j, 1\le j \le 1.5 \times 10^6\}$, where each $x_j$ is a ratio of the form (empirical second moment/predicted second moment) obtained using a block of $10^4$ zeros from $A23 \cup B23$. 
The sequence $\{x_j, 1\le j \le 1.5 \times 10^6\}$ was then sorted in descending order. 
That resulted in another sequence  $\{y_j, 1\le j\le 1.5 \times 10^6\}$; so, $y_1$ corresponds to the largest contribution among all blocks. The entries corresponding to $2k=2$ were then calculated as

\begin{equation}
100\times \frac{\sum_{j=1}^ n y_j}{\sum_{j=1}^{1.5 \times 10^6} y_j}\,, \qquad n=1,2,3,4,5\,. \
\end{equation}

\begin{table}[ht]
\footnotesize
\caption{\footnotesize For each $k$, samples of 10,000 zeros each are sorted according to their contribution to the $2k^{th}$ moment in decreasing order. For each $k$, the column lists the percentage contribution of the first $n$ sorted samples to the sum of all samples (1.5 million samples in total).}\label{na8} 
\begin{tabular}{|c|c|c|c|c|c|c|}
\hline
$n$ &\multicolumn{6}{|c|}{$A23 \cup B23$: 1.5M samples 10,000} \\
\cline{2-7}
 & $2k=2$ & $2k=4$ & $2k=6$ & $2k=8$ & $2k=10$ & $2k=12$ \\
\hline
1 & 0.0 & 0.1 & 0.8 & 4.0 & 10.0 & 17.2\\ 
2 & 0.0 & 0.1 & 1.6 & 7.6 & 18.9 & 32.4\\ 
3 & 0.0 & 0.1 & 2.1 & 9.9 & 24.2 & 40.6\\ 
4 & 0.0 & 0.2 & 2.6 & 11.8 & 28.0 & 46.1\\ 
5 & 0.0 & 0.2 & 3.0 & 13.4 & 31.4 & 50.8\\ 
\hline
\end{tabular} 
\end{table}

From Table~\ref{na8},  we see that near zero number $10^{23}$, more than 50\% of the value of the $12^{th}$ moment is determined by the 5 largest contributing samples out of a total of 1.5 million samples.  In comparison, the contribution of the analogous samples in the case of the $2^{nd}$ moment is negligible. This statistic illustrates the difficulty of sufficiently controlling high moments in experiments. 

\begin{figure}[ht]
\caption{
\footnotesize
Largest block contributions sorted in descending order and normalized according to (\ref{eq:largestcont})  using a sample of 1.5 million blocks from near $T=10^{22}$ (each block is $10^4$ zeros).}\label{powerlaw}
\includegraphics[scale=.75]{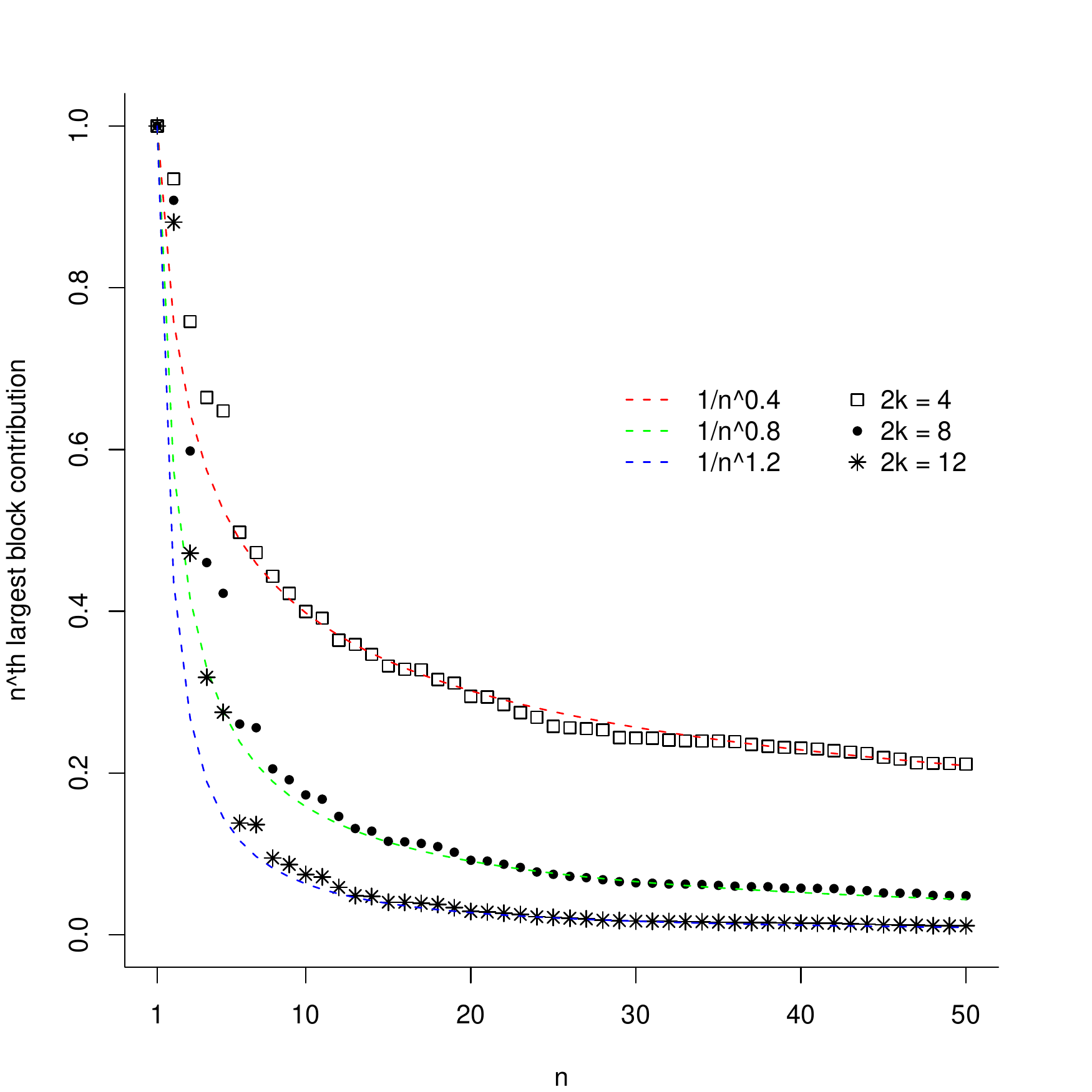}
\end{figure}

To further understand extreme values of the moments, we consider the function

\begin{equation} \label{eq:largestcont}
f(n):=f_k(n)=\frac{y_n}{y_1}\,, \qquad n \ge 1\,,
\end{equation}

\noindent 
which records the $n^{th}$ largest block contribution as a fraction of the largest block contribution, where each block spans $10^4$ zeros. It is worth mentioning that the largest contribution to the $12^{th}$ moment among such blocks found in our computations is 148,569 times the expectation according to Table~\ref{exp1}. Figure~\ref{powerlaw} is a graph of $f_k(n)$ for $1\le n\le 50$, and $2k=4$, 8, or 12. The figure is based on the full set of 15 billion zeros near $T=10^{22}$ (sets $A23$ and $B23$). 

Apart from a small initial segment, the lines in Figure~\ref{powerlaw} parallel rather closely the 
power law  $p(n):= p_k(n)= n^{-k/5}$. 
Moreover, we obtain a similar graph if we instead use the set of $10^9$ zeros 
near $T=10^{15}$ available to us (i.e. set $Z16$). 
Figure~\ref{powerlaw} would arise if the $n^{th}$ largest value of $|\zeta(1/2+it)|$ in the interval of integration
were about $n^{-1/10}$ times the largest value.  The log-normal distribution law suggests
that asymptotically this should hold with the 1/10 replaced by 

\begin{equation}
\frac{1}{2}\sqrt{\frac{\log \log T}{\log M}}\,.
\end{equation}

\noindent
where $M$ is the number of blocks. 

%



As evidenced by the data so far, the empirical moments vary substantially, even when calculated using long blocks of zeros. So, it might be interesting to consider

\begin{equation}\label{eq:apper21}
\log\left(\frac{\int_{ B_n}|\zeta(1/2+it)|^{2k}\,dt}{\int_{B_n} P_k\left(\log \frac{t}{2\pi}\right)\,dt}\right)\,. 
\end{equation}

\noindent
In other words, we consider $\log$(empirical moment/predicted moment), which is equal to 

\begin{equation}\label{eq:apper31}
\log\left(1+\frac{E_k(\beta_n)-E_k(\alpha_n)}{\int_{\alpha_n}^{\beta_n} P_k\left(\log \frac{t}{2\pi}\right)\,dt}\right)\,, 
\end{equation}

\noindent
where $\alpha_n$ and $\beta_n$ denote the ordinates of the first and last zeros of block $B_n$ respectively. The distribution of this quantity could shed some light on the remainder term in the full moment conjecture. Table~\ref{na10} contains  statistical summaries for (\ref{eq:apper31}) using several zero sets.  The numbers do not vary much from sample to sample. 
In particular, the various summary statistics change relatively little across the lower tables, 
which come from approximately the same height. 
Also, Table~\ref{na11} suggests if we sample (\ref{eq:apper31}) over many blocks $B_n$, 
 then normalize the samples to have mean 0 and variance 1, 
 the moments of the normalized samples generally start decreasing when $2k \ge 6$. 
It is not clear why there is such a trend, or whether it will persist. 


\begin{table}[ht]
\footnotesize
\caption{\footnotesize $\log$(empirical moment/predicted moment), where predicted moment is given in Table~\ref{exp1}. The columns Min, Max, Mean, and SD refer to the min log ratio, max log ratio, mean log ratio, and the standard deviation of the logs of ratios about their empirical mean. The set from $B23$ corresponds to $b2$ in Table~\ref{allsamples}.}\label{na10}
\begin{tabular}{|c||c|c|c|c||c|c|c|c||}
\hline
$2k$ &\multicolumn{4}{|c||}{$Z16$: 10,000 samples 100,000} & \multicolumn{4}{|c||}{$O20$: 10,000 samples 100,0000} \\
\cline{2-9}
 & Min  & Max & Mean & SD&  Min & Max & Mean & SD \\ 
\hline
2 & -0.048 & 0.069 & 0.000 & \textbf{0.013}& -0.066 & 0.155 & 0.000 & \textbf{0.021} \\ 
4 & -0.598 & 2.054 & -0.031 & \textbf{0.230}& -0.757 & 2.938 & -0.076 & \textbf{0.345} \\ 
6 & -1.871 & 4.400 & -0.326 & \textbf{0.671}& -2.468 & 5.298 & -0.657 & \textbf{0.886} \\ 
8 & -3.644 & 5.925 & -1.046 & \textbf{1.126}& -4.943 & 6.599 & -1.899 & \textbf{1.411} \\ 
10 & -5.749 & 6.917 & -2.134 & \textbf{1.560}& -7.931 & 7.232 & -3.666 & \textbf{1.907} \\ 
12 & -8.142 & 7.552 & -3.508 & \textbf{1.976}& -11.30 & 7.410 & -5.826 & \textbf{2.385} \\ 
\hline
\multicolumn{9}{c}{}\\
\hline
$2k$ &\multicolumn{4}{|c||}{$B23$: 10,000 samples 100,000} & \multicolumn{4}{|c||}{$A23$, $B23$: 150,000 samples 100,000} \\
\cline{2-9}
 & Min  & Max & Mean & SD&  Min & Max & Mean & SD \\ 
\hline
2 & -0.083 & 0.291 & 0.000 & \textbf{0.028} & -0.092 & 0.360 & 0.000 & \textbf{0.028} \\
4 &-0.949 & 3.938 & -0.121 & \textbf{0.423} & -1.067 & 4.343 & -0.120 & \textbf{0.418} \\
6 &-3.005 & 6.506 & -0.937 & \textbf{1.016} & -3.258 & 7.098 & -0.933 & \textbf{1.008} \\
8 &-5.977 & 7.847 & -2.580 & \textbf{1.578} & -6.335 & 8.623 & -2.573 & \textbf{1.568} \\
10 &-9.586 & 8.419 & -4.860 & \textbf{2.109} & -10.030 & 9.379 & -4.849 & \textbf{2.096} \\
12 &-13.82 & 8.462 & -7.613 & \textbf{2.623} & -14.330 & 9.606 & -7.597 & \textbf{2.607} \\ 
\hline
\end{tabular}
\end{table}



%
%


\begin{table}[ht]
\footnotesize
\caption{\footnotesize Moments of $\log$(empirical moment/predicted moment) normalized to have mean 0 and variance 1 (predicted moment given in Table~\ref{exp1}). The columns $p=3$, $p=4$, ...etc refer to the third moment, fourth moment,...etc, of the quantity $\log$(empirical moment/ predicted moment).}\label{na11}
\begin{tabular}{|c|c|c|c|c|c|c|}
\hline
$2k$ & \multicolumn{6}{|c|}{$Z16$: 10,000 samples 100,000}\\
\cline{2-7}
 & p=3 & p=4 & p=5 & p=6 & p=7 & p=8\\ 
\hline
2 & 0.271 & 3.338 & 3.490 & 22.64 & 51.80 & 264.8\\ 
4 & 1.553 & 8.237 & 39.06 & 238.3 & 1581 & 11380\\ 
6 & 1.306 & 6.020 & 21.40 & 99.97 & 489.3 & 2610\\ 
8 & 1.083 & 4.916 & 14.52 & 60.60 & 253.2 & 1172\\ 
10 & 0.954 & 4.424 & 11.68 & 46.68 & 178.9 & 774.5\\ 
12 & 0.875 & 4.170 & 10.23 & 40.21 & 146.1 & 611.1\\ 
\hline
\multicolumn{7}{c}{}\\
\hline
$2k$ &\multicolumn{6}{|c|}{$A23 \cup B23$: 150,000 samples 100,000} \\
\cline{2-7}
 & $p = 3$ & $p = 4$ & $p = 5$ & $p = 6$ & $p = 7$ & $p = 8$\\ 
\hline
2 & 0.862 & 5.855 & 25.53 & 200.2 & 1783 & 18240\\ 
4 & 1.737 & 8.773 & 42.40 & 255.8 & 1709 & 12560\\ 
6 & 1.286 & 5.747 & 19.83 & 89.95 & 430.8 & 2270\\ 
8 & 1.075 & 4.826 & 14.21 & 59.10 & 248.7 & 1167\\ 
10 & 0.966 & 4.445 & 11.99 & 48.34 & 190.6 & 850.4\\ 
12 & 0.902 & 4.251 & 10.84 & 43.22 & 164.0 & 712.7\\ 
\hline
\end{tabular}
\end{table}

Finally, we consider the correlations of the $2k^{th}$ moment. To quantify such correlations, we computed the autocovariances $c_m$ of the $2k^{th}$ moment. These are defined as

\begin{equation}\label{eq:autoco}
c_m := \frac{1}{R} \sum_{r=1}^R (x_{r+m} - \bar{x})(x_r - \bar{x})\,,
\end{equation} 

\noindent
where

\begin{eqnarray}
x_r := \int_{B_r} |\zeta(1/2+it)|^{2k}\,dt\,,\qquad \bar{x} =\frac{1}{R} \sum_{r=1}^R x_r\,,
\end{eqnarray}

\begin{figure}[ht]
\footnotesize
\caption{Correlations of the second moment.}\label{blockcorr}
\includegraphics[scale=.7]{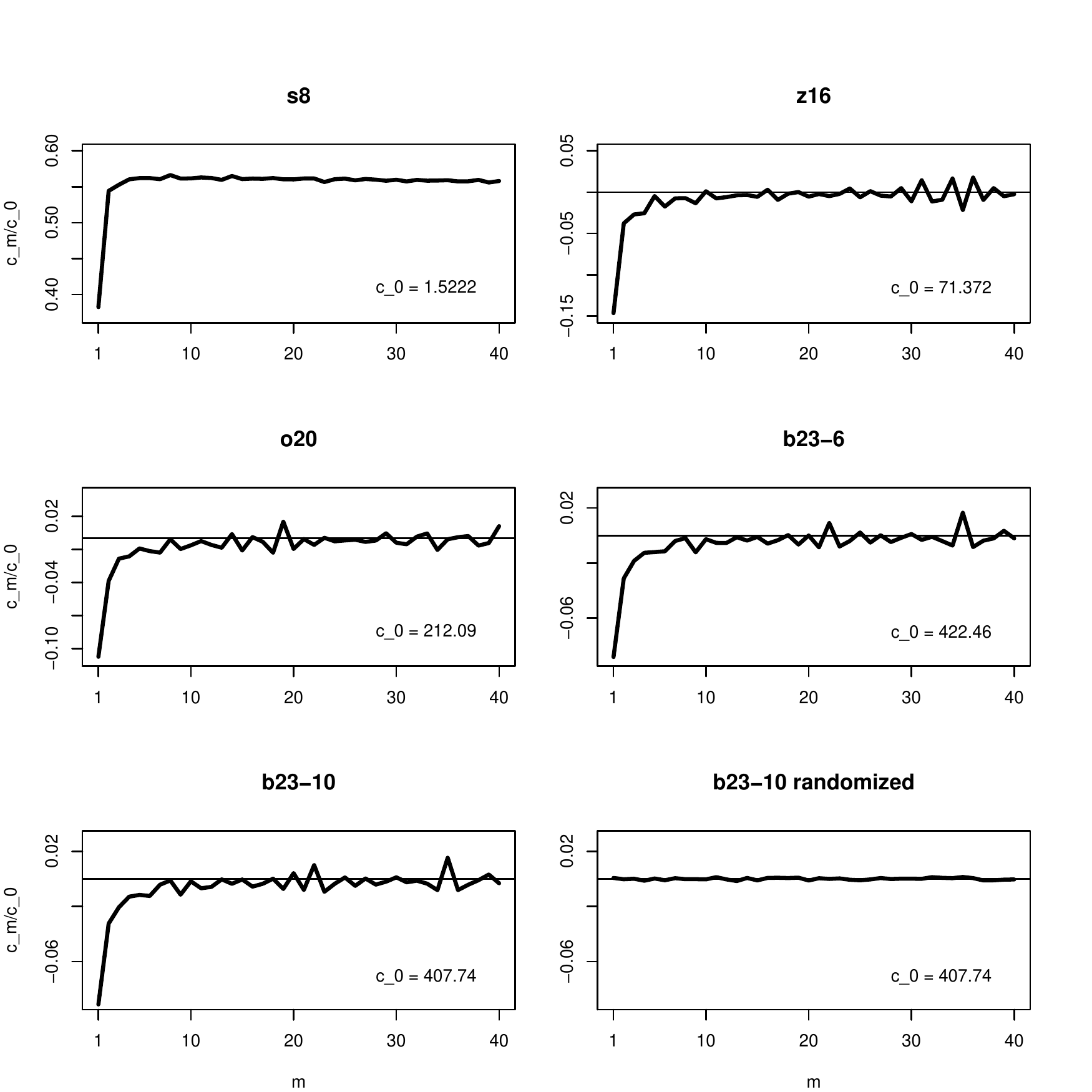}
\end{figure}

\noindent
We used  a total of $10^6$ consecutive blocks $B_r$ consisting of $10^3$ zeros each 
(so in definition (\ref{eq:autoco}), we have $R=10^6$). Figure~\ref{blockcorr} contains graphs of $c_m/c_0$ for $m=1,2,\ldots,40$, and $2k=2$, near heights $T=10^7$ (set s8), $T=10^{15}$ (set z16), $T=10^{19}$ (set o20), and $T=10^{22}$ (sets b23-6 and b23-10). 

Figure~\ref{blockcorr} has several interesting features. Although the autocovariances are overwhelmingly positive at low heights (set s8), they become mostly negative at large heights (e.g. set b23-10). Additionally, the amount of variation in the autocovariances $c_m$ is substantially more for the first few $m$ at all heights considered.

We plotted $c_m/c_0$ for two different sets near $T=10^{22}$ (sets b23-6 and b23-10). 
The two plots look almost identical, which suggests the autocovariances are significant. 
To better quantify this, we plotted the autocovariances of set b23-10 with the blocks 
randomly permuted (set ``b23-10 randomized''). 
Quite visibly, the autocovariances of the randomized set are much smaller than those of 
the original ordered set.  This points to some long range correlations in the values of the second moment. 

As for higher moments, the autocovariances do not appear to be significant. When $2k=4$, they are hardly distinguishable from those of a randomized set, and when $2k=6, 8, 10, 12$, they become completely indistinguishable.

The apparent significance of the autocovariances in the case of the second moment prompted us to consider the following shifted fourth moment:

\begin{equation}\label{eq:malpha}
M(T, H; \alpha) := \frac{1}{H} \int_T^{T+H} |\zeta(1/2+it)|^2\,|\zeta(1/2+it +i\alpha)|^2\,dt\,.
\end{equation} 

\noindent
If the values of $|\zeta(1/2+it)|$ and $|\zeta(1/2+it+i\alpha)|$ are statistically independent of each other, which is plausible for large $\alpha$, then one expects the integral in (\ref{eq:malpha}) to split, so it is approximated by

\begin{equation} \label{eq:malpha1}
\begin{split}
M(T,H; \alpha)=&\frac{1}{H}\int_T^{T+H} |\zeta(1/2+  it)|^2\, |\zeta(1/2+it+i\alpha)|^2\,dt \approx \\
&\left( \frac{1}{H}\int_T^{T+H} |\zeta(1/2+it)|^2\,dt\right) \,\left( \frac{1}{H}\int_T^{T+H} |\zeta(1/2+it+i\alpha)|^2\,dt\right)\,,
\end{split}
\end{equation}

\begin{figure}[ht]
\footnotesize
\caption{Shifted fourth moment $M(T,H; \alpha)/M(T,H; 0)$, black dots (drawn for $\alpha$ a multiple of $0.03$), versus the kernel $K(T; \alpha)$ defined in (\ref{eq:kernel}), dashed line. $T\approx 10^{22}$, $H\approx 6.5 \times 10^5$. }\label{smg1}
\includegraphics[scale=.7]{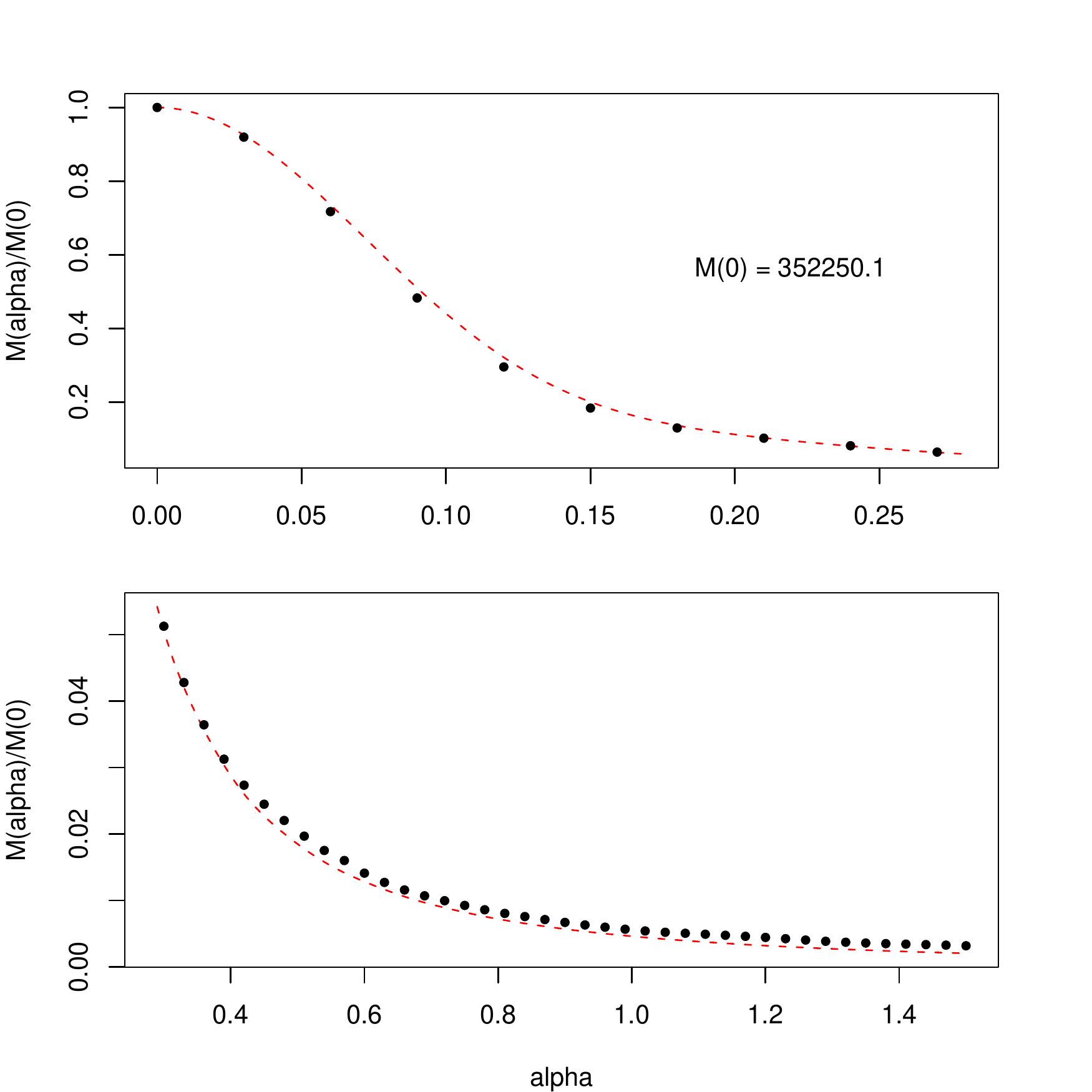}
\end{figure}
\begin{figure}[ht]
\footnotesize
\caption{Shifted fourth moment $M(T,H; \alpha)/M(T,H; 0)$, black dots (drawn for $\alpha$ a multiple of $0.5$), versus the kernel $K(T; \alpha)$ defined in (\ref{eq:kernel}), dashed line. $T\approx 10^{22}$, $H\approx 6.5 \times 10^5$. }\label{smg2}
\includegraphics[scale=.7]{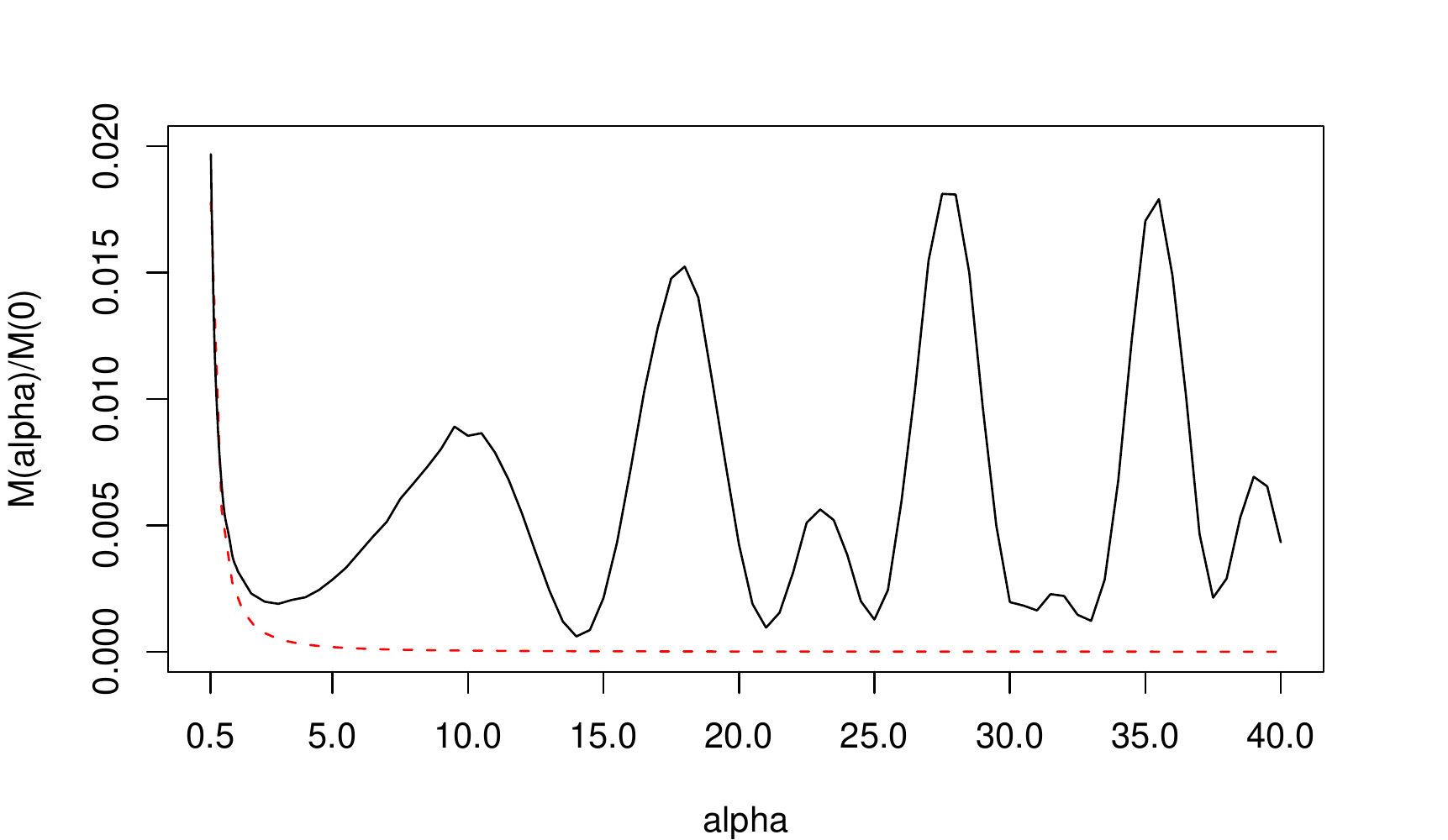}
\end{figure}

 
\noindent
and the latter is directly related to the autocovariances. K\"osters~\cite{Ko} and Chandee~\cite{Ch} investigated the function $M(T, T; \alpha)$ as well as other more general shifted moments. K\"osters' work immediately implies if $\alpha \log T = O(1)$  (as $T\to \infty$), and if

\begin{equation}\label{eq:kernel}
K(T; \alpha):=\frac{12}{(\alpha\,\log T)^2}\left(1-\frac{4 \sin^2 (\alpha\,\log T/2)}{(\alpha\,\log T)^2}\right)\,,
\end{equation}

\noindent
then

\begin{equation}\label{eq:malpha2}
\frac{M(T, T; \alpha)}{M(T, T; 0)}  \sim K(T; \alpha)\,,
\end{equation} 

\noindent
where right side is defined for $\alpha=0$ by taking a limit. It is plausible that relation (\ref{eq:malpha2}) would continue to hold if the left side is replaced by 

\begin{equation} \label{eq:malpha3}
\frac{M(T, H; \alpha)}{M(T, H; 0)} = \frac{\int_T^{T+H} |\zeta(1/2+it)|^2\,|\zeta(1/2+it + i\alpha)|^2\,dt}{\int_T^{T+H} |\zeta(1/2+it)|^4\,dt}\,,
\end{equation}

\noindent
when $H$ is much smaller than $T$ but not too small. 
This is supported to some extent by Figure~\ref{smg1} since it shows a reasonable agreement 
between (\ref{eq:malpha3}) and the kernel (\ref{eq:kernel}) for $0 \le \alpha \le 1.5$, 
$T\approx 10^{22}$, and $H\approx 6.5 \times 10^5$.  However, for $\alpha \ge 0.4$, Figure~\ref{smg1} points to  
 deviations from the kernel model.
 And Figure~\ref{smg2} shows that, starting around $\alpha =1.5$, there is a clear 
departure from the kernel model, and the onset of asymptotic oscillations (of amplitude $\le 0.018$). 
Furthermore, the oscillations remain significant for large $\alpha$ (large relative to $\log T$); 
for example, we computed $M(T, H; 100) = 0.011447$ and $M(T, H; 1000) = 0.004035$. 
Presumably, these long-range oscillations are induced by prime sums present in the lower order terms of relation (\ref{eq:malpha2}). 





\section{Numerical methods}

We calculated the moments of $|\zeta(1/2+it)|$ by evaluating integrals similar 
to (\ref{eq:sndconj}) near $T=10^{22}$ and at lower heights. 
The computational method involved two main choices,
 a method to approximate point-wise values of $|\zeta(1/2+it)|$ for large $t$, and  
 a suitable integration technique to handle integrals like (\ref{eq:sndconj}). 

\subsection{Point-wise approximations}

One method to numerically evaluate \mbox{$\zeta(1/2+it)$} is based on 
Euler-Maclaurin summation. This method is derived from the summation by 
parts formula, and it can evaluate the zeta function on the critical 
line with great accuracy.  But it requires $t^{1+o_{\kappa}(1)}$ elementary operations on numbers of $O_{\kappa}(\log t)$ bits  to compute $\zeta(1/2+it)$ to within $\pm t^{-\kappa}$ for a 
single value of $t$ even if we do not demand high accuracy. 
This is prohibitive when $t$ is near $10^{22}$, as in our experiments. 
Note that asymptotic constants are taken as $t\to \infty$, and the notations 
$o_{\kappa}()$ and $O_{\kappa}()$ mean asymptotic constants depend on the parameter $\kappa$ only. 

A much more efficient method to numerically approximate $\zeta(1/2+it)$ is the Riemann-Siegel formula, which was invented by 
Riemann and rediscovered by Siegel during the latter's study of Riemann's notes. 
The Riemann-Siegel formula can be used to compute $\zeta(1/2+it)$ to within $\pm t^{-\kappa}$ for 
large $t$ using $t^{1/2+o_{\kappa}(1)}$ operations on numbers of $O_{\kappa}(\log t)$ bits. 


The running times stated above are for computing single values of the zeta
function.
For computations that involve many evaluations of $\zeta(1/2+it)$ in a restricted interval, the Odlyzko-Sch\"onhage~\cite{OS} algorithm (OS) offers 
significant savings.

\begin{thm} [The Odlyzko-Sch\"onhage algorithm (OS)] \label{thm:os}
Given $c>0$, $\epsilon >0$, and $a \in [0,1/2]$ there exists an effectively computable 
$c_1=c_1(c,a,\epsilon)>0$ so that one can compute $Z(t)$ for any value of $t>0$ in the 
range $T<t<T+T^a$ with accuracy $T^{-c}$ using $\le c_1t^{\epsilon}$ operations on 
numbers of $\le c_1\log t$ bits provided a precomputation involving 
$\le c_1T^{1/2+\epsilon}$ operations and $\le c_1T^{a+\epsilon}$ bits of storage is done beforehand. 
\end{thm}

For example, choosing $a=1/2$ in Theorem~\ref{thm:os}, one can compute 
$Z(t)$ at $n\approx T^{1/2}$ points in the range $(T,T+\sqrt{T})$ using 
$n^{1+\epsilon}$ operations, which is nearly optimal. By comparison, the Riemann-Siegel formula requires $n^{2+\epsilon}$ operations to achieve the same task, and the Euler-Maclaurin formula requires $n^{3+\epsilon}$ operations.

Our approach to computing values of the zeta function was to use the
results of the precomputations of the OS algorithm
that had been carried out by the second author at AT\&T Labs - Research.

Although computations using the OS method, as well as the other methods mentioned previously,
can be made completely rigorous, provided one uses multi-precision arithmetic,
 this is almost never done in practice because of the time cost.
Instead, the validity of the
final results depends on the assumption that there is quasi random cancellation of round-off
errors, plus some additional checks.  For a discussion, see \cite{HO,O1,O2}.

In addition to OS, we used another method to compute $|\zeta(1/2+it)|$ that 
 appears to be valid over short ranges and is numerically analyzed in detail is Section 6. 
 To describe briefly, though, suppose in some zero interval $I_n:=(\gamma_n,\gamma_{n+1})$ one knows the 
value of $|\zeta(1/2+i(\gamma_n+\gamma_{n+1})/2)$ as well as the location of the $m$ zeros below and
$m$ zeros above $\gamma_n$. 
Then one can try to approximate $|\zeta(1/2+it)|$ for  $t\in I_n$ by a polynomial that assumes the value 
$|\zeta(1/2+i((\gamma_n+\gamma_{n+1})/2)$ at $(\gamma_n+\gamma_{n+1})/2$ and vanishes at those $2m$ neighboring zeros. 
Let us denote this approximation by HP since it is essentially a truncated Hadamard product. 
Numerical experiments, which are discussed in Section 5, suggest the quality of HP 
improves linearly with $m$, which is the number of zeros used (the improvement is in the sense of $L^{\infty}$-error). 

In general, HP approximations are not very accurate. 
For example, with $m=2^9$ and for $t\approx 10^{22}$, HP has an $L^{\infty}$ error 
of about $6\times 10^{-2}$.  On the other hand, HP is faster than OS. If used 
judiciously (see Section 5.3 for a detailed explanation), it can reduce the running 
time of our computations while not affecting the expected accuracy of the moment data. 

\subsection{Integration Methods}

Our goal is to calculate the 2nd to 12th even moments of $|\zeta(1/2+it)|$.\footnote{Our computations were slightly more comprehensive because they included the $-0.5$ to $12$ moments in increments of 1/2. But only the even moments are discussed in this paper.} 
To control for oscillations, each integration interval has two consecutive zeros of 
$Z(t)$ as its endpoints; i.e., is of the form $I_n=(\gamma_n,\gamma_{n+1})$. 

We experimented with several integration methods. We settled on Romberg integration as our choice, 
due to its simplicity, fast convergence, and ability to naturally provide a  
posteriori error estimates. 

\subsection{Numerical Implementation}

We first discuss implementation details in the vicinity of zero number $10^{23}$ 
(i.e. sets $A23$ and $B23$). Numerical tests on overlapping intervals showed the distribution of the $L^{\infty}$-error of the OS approximation of $|\zeta(1/2+it)|$ is normally distributed with mean $\approx 0$ and standard deviation $\approx 5\times 10^{-9}$, so the OS approximation is typically accurate to within $\pm 5\times 10^{-9}$. In order to determine the integration intervals, we used the zero sets previously compiled by 
the second author. Numerical tests on overlapping zero sets showed the $L^{\infty}$-error in the location of zeros 
hovered around $5\times 10^{-8}$. 

Our initial goal was to obtain enough accuracy to enable comparisons 
between the empirical moments and prediction (\ref{eq:conjmoments}).  For each zero 
interval, we attempted to compute the $2k^{th}$ moment so that

\begin{equation}
\footnotesize
\begin{split}
\textrm{Absolute posteriori error} < & 10^{-3}\times  \textrm{Average zero gap at height $T$} \ \\
&\times \textrm{Expected $2k^{th}$ moment according to (\ref{eq:conjmoments})}\,. 
\end{split}
\end{equation}

\noindent
Note the average zero gap near height $T=10^{22}$ is about $0.128$. 
We aggregated the moment data for each 1,000 zeros.  
Thus, we expect the aggregate error in our computed values of the $2k^{th}$ moment divided by the length of integration interval to satisfy

\begin{equation}
\footnotesize
\textrm{Aggregate posteriori error}<10^{-3}\times \textrm{Expected $2k^{th}$ moment according to (\ref{eq:conjmoments})}\,. 
\end{equation}

\noindent
In particular, upon dividing the empirical moment by prediction (\ref{eq:conjmoments}), the ratio will typically have 3 significant decimal digits. 
For the values of $k$ we are considering though, the leading term predictions (\ref{eq:conjmoments}) are less in size than the full moment predictions  (\ref{eq:shortmoment}). 
Since the latter are conjectured to be more accurate, then upon dividing empirical moments by the full moment prediction (\ref{eq:shortmoment}), the ratio should have more digits of accuracy (especially for high moments). That is, we expect the ratio (empirical moment (\ref{eq:sndconj})/predicted moment (\ref{eq:shortmoment})) to be correct to within

\begin{equation}
\footnotesize
\pm 10^{-3}\times \frac{\textrm{Expected $2k^{th}$ moment according to (\ref{eq:conjmoments})}}{\textrm{Expected $2k^{th}$ moment according to (\ref{eq:shortmoment})}}\,. 
\end{equation}

Once we decided on this accuracy standard, the following were chosen accordingly:  the 
number of zeros to use in HP, and a criterion to determine whether HP or OS is more appropriate. 
After some numerical tests, which are described in Section 4, we decided to fix the 
number of zeros supplied to HP at 1000 zeros, or 500 zeros on each side of each interval. 
With this choice, HP is about 5 times faster than OS, it has an $L^{\infty}$-error of about $6\times 10^{-2}$ and a typical error (i.e. square root of $L^2$-error) of about $10^{-4}$. 

Given the unimodal shape of $|\zeta(1/2+it)|$ on each integration interval $I_n=(\gamma_n,\gamma_{n+1})$, it is reasonable to try to approximate $\max_{t\in I_n} |\zeta(1/2+it)|$ by $C:=|\zeta(1/2+it_0)|$, where $t_0$ is the midpoint of $I_n$. So, a rough upper bound for the contribution of $I_n$ to the $2k^{th}$ moment is $C^{2k}$. In particular, if $C$ is small enough, then HP can be safely used; i.e. without affecting the accuracy standard. Numerical experiments suggested that given our choices so far, we could set the $C$-threshold for using HP at $C\le 7.0$. Note that the contribution of such intervals is, in any case, mostly negligible for the $4^{th}$ and higher moments at all heights considered.

To test the accuracy and efficiency of our implementation, we set it to work in regions where $|\zeta(1/2+it)|$ assumes large values. Based on this, we chose an upper bound of 2048 iterations for Romberg integration.  The integration programs used  double arithmetic. 
This allows for about $15$ to $16$ significant digits, which is more than sufficient for our purposes. 

In sum, with our choices of the parameters, the program used HP to
integrate on about 88\% of the intervals, the posteriori error typically surpassed 
the accuracy standard, and the program consumed an average of 12 OS function evaluations 
per interval\footnote{We counted an HP evaluation as 1/5 OS evaluation.}. 
For lower sets $O20$ and $Z16$, the parameter choices were decided similarly.  
As for set $S8$, which consists of the first $10^8$ zeta zeros, we relied exclusively on the Riemann-Siegel formula to compute the moments. Lastly, as a check of the validity of the  integration programs, we successfully reproduced moment data computed by Michael Rubinstein for $t\approx 10^6$. We implemented the code in FORTRAN 90, and ran it onSGI machines at the Minnesota Supercomputing Institute.

\section{Numerical tests of local models of the zeta function}

In addition to HP, which is the polynomial approximation discussed in Section 5, there is another attractive formula to approximate $\zeta(1/2+it)$ over a zero interval $I_n=(\gamma_n,\gamma_{n+1})$ which is due to Gonek, Hughes, and Keating~\cite{GHK}.
The method expresses $\zeta(1/2+it)$ as the product of two parts: one that resembles a truncated Euler product, and another that resembles a truncated Hadamard product. 
For this reason, we will abbreviate the~\cite{GHK} formula as EHP (Euler-Hadamard product). 
The approximation EHP does not require the value of $\zeta(1/2+it)$ at the midpoint of the 
zero interval $I_n$. 
Instead, it incorporates the contribution of the primes in a natural way. 
The approximation arises from a smoothed approximate functional equation for the logarithmic derivative of zeta. 

EHP requires the Euler and Hadamard products to be both suitably smoothed. 
We implemented smoothing for the Euler product, but not for the Hadamard product. 
It was pointed out to the authors, however, that it is difficult to eliminate smoothing 
from the Hadamard product in EHP with realistic error bounds~\cite{H}. 
This suggests smoothing the Hadamard product is critical to the accuracy. 
Nevertheless, implementing such smoothing would cause EHP to be too slow in practice. 
Thus, in considering whether to use EHP in our moment calculations, we only tested it 
with an unsmoothed Hadamard product, which is also how it was briefly tested in~\cite{GHK}. 
We found the accuracy of EHP to be comparable to HP, but the latter was significantly 
faster due to its simplicity (the convergence of HP is linear in the number of zeros used).  

\subsection{Definitions}

Gonek et al. have obtained the following formula for the Riemann zeta function:

\begin{equation} \label{eq:one}
\zeta(s)=P_X(s)Z_X(s)\left[1+O\left(\frac{X^{K+2}}{(|s|\log X)^K}\right)\right]\,,
\end{equation}

\noindent
where $s=\sigma+it$, $\sigma\ge 0$, $|t|\ge 2$, $X\ge2$, K any fixed positive integer, and

\begin{equation} \label{eq:two}
P_X(s)=\exp\bigg(\sum_{n\le X}\frac{\Lambda(n)}{n^s\log n}v(e ^{\log n/\log X})\bigg)\,,
\end{equation}

\begin{equation} \label{eq:three}
Z_X(s)=\exp\bigg(-\sum_{\rho}U((s-\rho)\log X)\bigg)\,,
\end{equation}

\noindent
where $u(x)$ is a nonnegative $C^\infty$ function of mass 1 with support on [$e^{1-1/X},e$], $v(t)=\int_t^\infty u(x)\,dx$, $U(s)=\int_0^{\infty}u(x)E_1(s\log x)\,dx$, $\rho$ denotes a nontrivial zeros of the zeta function, and $E_1(z)=\int_z^{\infty} e^{-w}/w\,dw$ is the exponential integral. The $O$-notation constants depend on $u$ and $K$ only.

A literal implementation of (\ref{eq:one}) is not efficient because it is expensive to compute $U(s)$. For this reason, we use an unsmoothed, truncated version of (\ref{eq:three}); in particular, we replace $u(x)$ in the definition of $U(s)$ by a delta function at $e$. We assume the Riemann hypothesis, and for $t \in (\gamma_n,\gamma_{n+1})$  we define

\begin{equation} \label{eq:four}
Z_X^{n,m}(1/2+it):=\exp\bigg(-\sum_{j=n-m+1}^{n+m}E_1\big(i(t-\gamma_j)\log X)\big)\bigg)\,.
\end{equation}

As for $P_X(1/2+i t)$, we did apply smoothing to it because  smoothing there is not computationally expensive. Specifically, we first define the following non-negative $C^{\infty}$ function:

\begin{displaymath}
f(x)=\left\{
\begin{array}{cr} 
e^{-1/x^2} & x>0 \\ 
0 & x<0
\end{array}
\right.,\qquad g(x)=\frac{f(x)f(1-x)}{\int_0^1f(x)f(1-x)\,dx}\,,
\end{displaymath} 

\noindent
then define the smoothing kernel

\begin{equation} \label{eq:six}
u(x)=X g(X\log(x/e) + 1)/x\,.
\end{equation}

\begin{figure}[ht]
\caption{\footnotesize The kernel $u(x)$ with $X=6$} \label{pppxxx}
\includegraphics[scale=0.8, angle=0, viewport=90 575 325 720,clip]{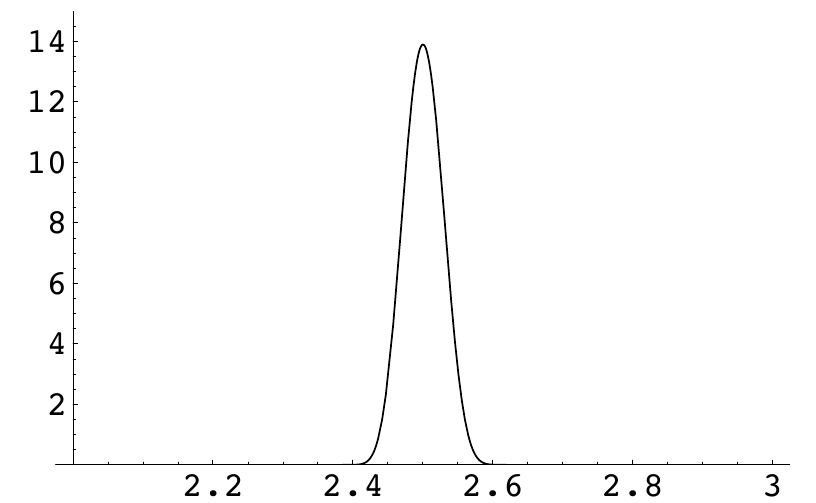}
\end{figure}

\noindent
Figur~\ref{pppxxx} is a plot of $u(x)$ when $X=6$. Put together then, for $t \in$ ($\gamma_n,\gamma_{n+1}$) we have

\begin{equation}\label{eq:ehpapp2}
\begin{split}
|\zeta(1/2+it)|&=:|P_X(1/2+ it)Z_X^{n,m}(1/2+it)| +R_{1,m}(t) \\
& \\
&=\exp\bigg(\displaystyle \sum_{ n\le X}\frac{\displaystyle \Lambda(n)\cos(t\log n)}{\displaystyle \sqrt{n}\log n}v(e ^{\log n/\log X})\bigg) \\
&\,\,\,\,\,\,\,\, \times \exp\bigg(\displaystyle \sum_{j=n-m+1}^{n+m}Ci\big(|t-\gamma_j|\log X\big)\bigg) +R_{1,m}(t)\,,
\end{split}
\end{equation}

\noindent
where $Ci(x)=\int_x^{\infty} \cos(t)/t\,dt$ is the cosine integral, and $R_{1,m}(t)$ is the remainder function. We remark the numerics showed little sensitivity to the exact choice of $u(x)$ (which determines $v(x)$ in (\ref{eq:ehpapp2}). 

Let us denote the approximation $|P_X(1/2+ it)Z_X^{n,m}(1/2+it)|$ in (\ref{eq:ehpapp2})  by EHP since it is an Euler-Hadamard product.  Note the approximation depends on two parameters: $X$ and $m$, which control the number of primes and zeros used in $P_X$ and $Z_X^{n,m}$ respectively. 

In addition to EHP, we considered the following approximation of $|\zeta(1/2+it)|$: let

\begin{displaymath} Q^{n,m}(t)=\displaystyle
\prod_{j=n-m+1}^{n+m}|t-\gamma_j|\,,
\end{displaymath}

\noindent
then, for $t \in$ ($\gamma_n,\gamma_{n+1}$), we define

\begin{equation} \label{eq:eight}
|\zeta(1/2+it)|= Q^{n,m}(t)\, \frac{|\zeta(1/2+i\eta_n)|}{Q^{n,m}(\eta_n)} + R_{2,m}(t)\,,
\end{equation}

\noindent
where $\eta_n=(\gamma_{n+1}+\gamma_n)/2$, and $R_{2,m}(t)$ is the remainder function. We denote the approximation $ Q^{n,m} (t)\,\frac{Z(1/2+i\eta_n)|}{Q^{n,m}(\eta_n)}$ by HP as it is a truncated Euler product. We expect it to be a good approximation of $|\zeta(1/2+it)|$ because locally, away from the pole at 1,  the zeta function may be treated as a polynomial with the non-trivial zeros as its roots. 

\subsection{The function $|P_X(1/2+it)|$}

\begin{figure}[htbp]
\caption{\footnotesize $|P_X(1/2+i(T+t))|$ for $X=6$ (solid), 50.92 (short dashes), and 1000 (long dashes). The interval $(T,T+5)$ covers about 40 zeros.\label{varpx}}
\includegraphics[scale=0.8, angle=0, viewport=93 575 325 720,clip]{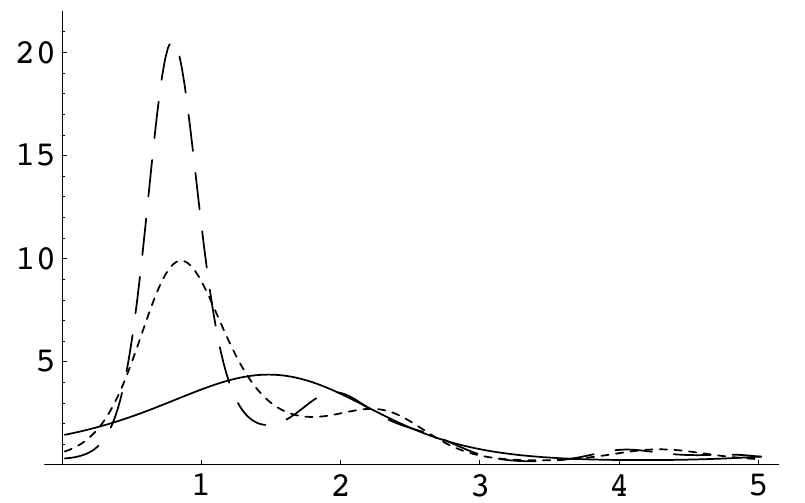}
\end{figure}

Formula (\ref{eq:one}) incorporates arithmetic contributions to the zeta function via $P_X$. 
We tested whether $P_X$ can be expected to have a tangible impact on numerical results by measuring its variation about its mean. Figure~\ref{varpx} is a plot of $|P_X(1/2+iT)|$ near 
 
\begin{equation} \label{eq:heighttt}
T=1.306643440879589721233593307594 \times 10^{22}\,, 
\end{equation}

\noindent
which is in the vicinity of the $10^{23}$-rd zero.  It shows $|P_X(1/2+it)|$ varies substantially. 
Furthermore, we calculated the variance of $|P_X(1/2+it)|$ in an interval of the form $(T,T+H)$; that is,  

\begin{equation}
\frac{1}{H}\int_0^H \left||P_X(1/2+i(T+t))|-\overline{P_X}\right|^2\,dt\,,\qquad \overline{P_X}:=\frac{1}{H}\int_0^H|P_X(1/2+i(T+t))|\,dt\,.
\end{equation}

\noindent
We chose $H=128$, which approximately the length of the stretch covered by 1,000 zeros near height $T$.  Table~\ref{pxt}  lists the variance for several values of $X$. It shows in particular  both the mean and the variance of $P_X$ increase with $X$. 

\begin{table}[ht] 
\footnotesize
\caption{\footnotesize Variation of $|P_X(1/2+i(T+t)|$ for $t\in (0,128)$. \label{pxt}}
\begin{tabular}{lccc}\hline
Value of $X$ & 6 & 50.92 &1000\\ 
Mean & 1.36 & 1.67 & 1.93 \\ 
Variance & 1.33 & 4.88 & 9.68\\ \hline
\end{tabular}
\end{table}

\subsection{Approximations of $|\zeta(1/2+it)|$}

\begin{figure}[htbp]
\caption{\footnotesize $|\zeta(1/2+i(T+t))|$ (solid) and EHP (dashed) with $X=6$ and $m=25$ over the interval covering zeros numbered 865 to 885 in the set $W$.\label{tathree}}
\includegraphics[scale=0.8, angle=0, viewport=93 575 325 720,clip]{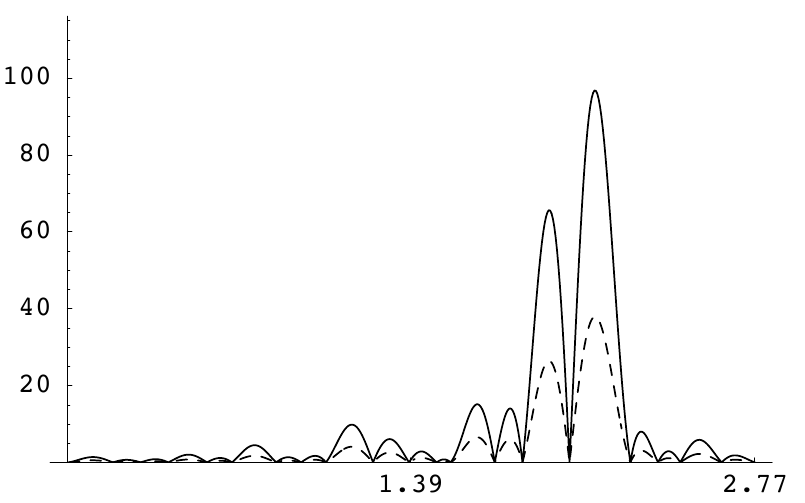}
\end{figure}

\begin{figure}[htbp]
\caption{\footnotesize $|\zeta(1/2+i(T+t))|$ (solid) and EHP (dashed) with $X=50.92$ and $m=16$ over the interval covering zeros numbered 865 to 885 in the set $W$.\label{ffive}}
\includegraphics[scale=0.8, angle=0, viewport=93 575 325 720,clip]{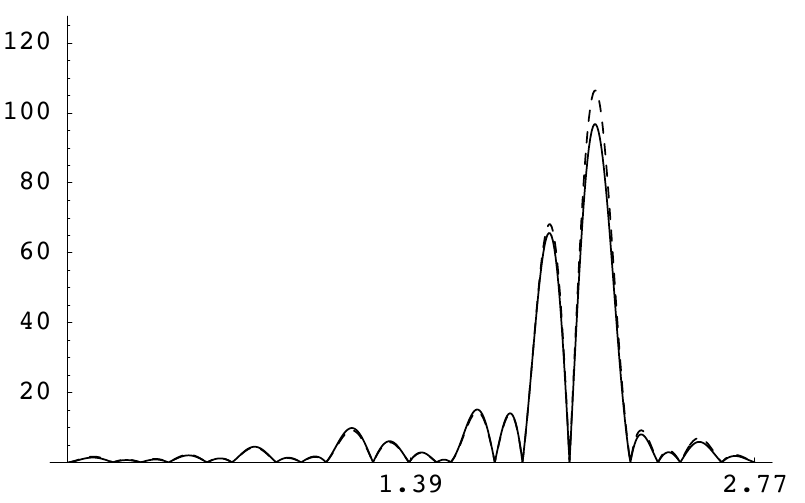}
\end{figure}

Our experiments relied on a set of 30,000 consecutive zeros near zero number $10^{23}$. We denote this set by $W$. The first zero in set $W$ is at the height $T$ specified in (\ref{eq:heighttt}). This is the first zero above Gram point number $10^{23} + 18,767,166,306$. Figures~\ref{tathree} and~\ref{ffive} are plots of the approximation EHP, which was defined in (\ref{eq:ehpapp2}), for several choices of $X$ and $m$ together with a plot of $|\zeta(1/2+it)|$  (calculated using the OS algorithm). The figures show that agreement is fair given the basic nature of the input supplied to EHP.

\subsection{Numerical experiments}
For each interval $(\gamma_n,\gamma_{n+1})$ in the set $W$, we calculated the values of $|\zeta(1/2+it)|$ at $l(n)$ equally spaced points $t_{n,d}$ inside $(\gamma_n,\gamma_{n+1})$, where $l(n)=2\lfloor 10\frac{\gamma_{n+1}-\gamma_n}{\Delta_n}+1\rfloor +1$, and $\Delta_n$ is the average spacing of the zeros near zero ordinate $\gamma_n$. Notice $l(n)$ is always odd, so $|\zeta(1/2+i(\gamma_{n+1}+\gamma_n)/2)|$, which is the value at the midpoint,  is always included. The values of $\zeta(1/2+it)$ were computed using the OS algorithm. We carried out 29 experiments, each involving 1,000 consecutive zeros;  that is, in experiment $v$, we used either HP or EHP to approximate $|\zeta(1/2+it)|$ over the interval 

\begin{equation}
[r_{v-1},r_v),\qquad r_v=\gamma_{N+500+1000v}, \qquad  v=1,\ldots,29\,, 
\end{equation}

\noindent
where $N\approx 10^{23}$ is the number of the first zero in set $W$.  In each experiment, we let $m$ (in $Z^{n,m}_X$ and $Q^{n,m}$) take on the values $2^r$ for $r=1,\ldots,8$, and we let $X$ take on the values 2, 6, 50.92, 1000, 2000, and 4000. Lastly, 
we let $Z^*(1/2+it)$ denote one of the approximations HP or EHP, and define the $L^{\infty}$-error in experiment $v$ by

\begin{equation}
\max_{\substack{d=1,\ldots,l(j) \\ j=r_{v-1},\ldots,r_v-1}} |\zeta(1/2+it_{j,d})-Z^*(1/2+it_{j,d})|\,. 
\end{equation}

\begin{figure}[htbp]
\caption{\footnotesize The $L^\infty$ error of HP (horizontal axes) vs.  EHP (vertical axes) with $m=256$ and $X=1000$ for all 29 experiments. The dashed line is $y=x$. \label{tafour}}
\includegraphics[scale=0.8, angle=0, viewport=93 575 325 720,clip]{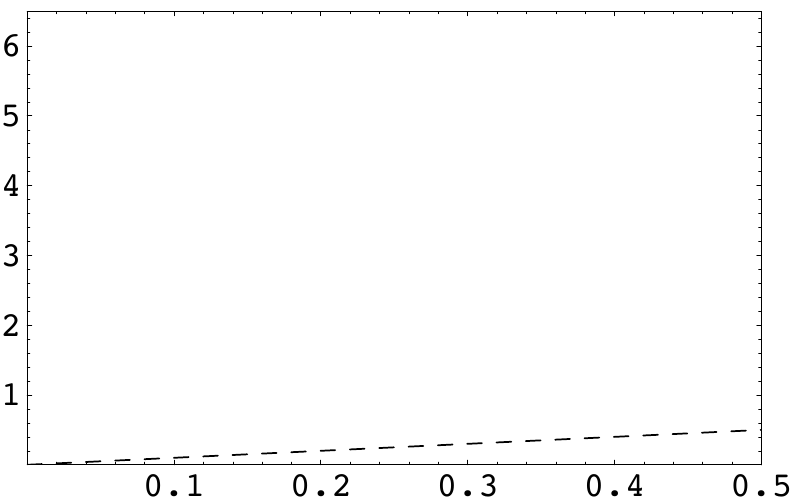}
\end{figure}

\noindent
Figure~\ref{tafour} is a scatter plot of the $L^{\infty}$-errors obtained with $m=64$ and $X=1000$. The line $y=x$ (dashed) is included to aid in comparison. The plot suggests that with the current choices of parameters $X$ and $m$, HP is generally a better approximation than EHP.

We are also interested in understanding the effect of modifying the parameter $X$, which essentially controls the importance of the arithmetic part $P_X$ relative to the polynomial part $Z_X$. Specifically, as $X$ decreases, $P_X$ becomes less important and we approach the HP model. To this end, consider Table~\ref{ttwo},\footnote{
For $X=50.92, 1000$ and 2000, we averaged the $L^\infty$ error with $m=256$ using the results from all of the 29 experiments. However, for $X=4000$, the average is based on data obtained from 15 experiments only.} 
 which indicates the $L^{\infty}$-errors reach a minimum when $X$ is in an intermediary region ($X\approx 1000$). But this minimum is 2.82, which is significantly larger than that of HP whose average is 0.13 with $m=256$.  Also, the standard deviation of the errors is roughly 1 for the values of $X$ listed in Table~\ref{ttwo}, whereas it is 0.07 for HP.  

\begin{table}[ht]
\footnotesize
\caption{\footnotesize The average $L^\infty$-error of EHP with $m=256$ over all 29 experiments.\label{ttwo}}
\begin{tabular}{lccccc}\hline
Value of $X$ &6& 50.92 & 1000 & 2000 & 4000 \\ 
Average $L^\infty$-error &6.84&3.81& 2.82 & 3.11 & 3.00\\ \hline
\end{tabular}
\end{table}

In a sense, the kind of comparison we have carried out so far may be ill-defined, because by construction HP is exact at $(\gamma_n+\gamma_{n+1})/2$, whereas EHP is not necessarily so. If EHP is normalized so that it is exact at the midpoint, the results improve significantly. Let us call this approximation ``normalized EHP.'' Table~\ref{tthree} contains the average $L^{\infty}$-errors of that approximation over all experiments. 
 The  $L^{\infty}$-errors for normalized EHP appear to reach a local minimum when $X\approx 6$ (i.e. the arithmetic part $P_X$ uses only the primes 2, 3, and 5). Also, table~\ref{tfour} contains the history of convergence of the $L^{\infty}$-errors as $m$ (which determines the number of zeros used) increases. 

\begin{table}[ht]
\footnotesize
\caption{\footnotesize The average $L^{\infty}$-error of normalized EHP with $m=256$ over all 29 experiments.\label{tthree}}
\begin{tabular}{lcccccc}\hline
Value of $X$ & 2 & $2.71828$ & 2.9 & 6 & 50.92 & 1000 \\ 
Average L$^\infty$-error &0.52& 0.16&0.30&0.11 &0.50 & 1.03 \\ \hline
\end{tabular}
\end{table}


\begin{table}[ht]
\footnotesize
\caption{\footnotesize History of Convergence of the $L^{\infty}$-error. Results based on experiment 1 (i.e. zeros numbered 500 to 1500 in set $W$).\label{tfour}}
\begin{tabular}{|c|c|c|c|} \hline
$m$ & HP & normalized EHP with $X=50$ & normalized EHP with $X=6$ \\ \hline
2   &   16.0327267  &     10.3965279                &      13.2324066              \\
4   &   8.6654952   &     7.3793465                &     4.9913900               \\
8   &    4.9345427  &         3.3031051             &  1.9168742                 \\
16   &   2.0833875  &   3.8302768                 &       2.7077135             \\
32   &   1.0230760   &        2.1464398             &     1.2456626               \\
64   &   0.5268745    &               1.6312318       &          0.5550487              \\
128   &  0.2601180   &      0.9890934                &   0.3347868                  \\
256   &   0.1308390   &     0.4652148                &     0.0544018    \\ \hline            
\end{tabular}
\end{table}

Finally, we considered the convergence rates of normalized EHP  for various values of $X$. The convergence rate is defined by

\begin{equation}
C_r=\frac{\log E_r-\log E_{r-1}}{\log (1/2)}, \qquad \textrm{where $E_r$ is the average $L^{\infty}$-error with $m=2^r$}\,. 
\end{equation}

\begin{figure}[ht]
\caption{\footnotesize The vertical axes is the $\log$ of the average $L^\infty$-errors over all 29 experiments of HP (solid),  EHP with $X=1000$ (short dashes), and normalized EHP with $X=6$ (long dashes). The horizontal axes is $\log m/\log 2$, where $m=2,4,8,\ldots,256$. The convergence rate is $\approx -slope\times\frac{1}{\log 2}$\label{fthirteen}}
\includegraphics[scale=1, angle=0, viewport=93 575 325 720,clip]{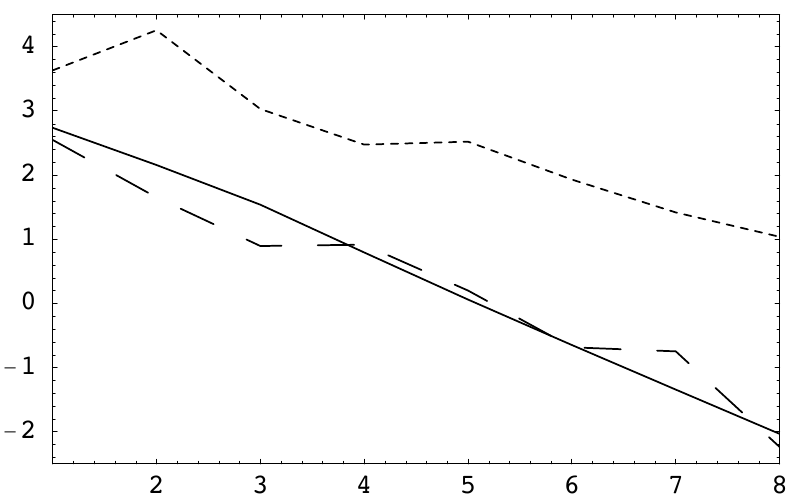}
\end{figure}

Table~\ref{tfive} indicates the convergence rate of normalized EHP is not smooth. 
Even orders of convergence fluctuate considerably, and sometimes are negative.
On the other hand, HP converges more smoothly at a linear rate. Figure~\ref{fthirteen} is a visual representation of this. We note the amount of time required in the EHP model was significantly more than in the HP model. Most of the additional time went into evaluating the cosine Integral . 

\begin{table}[ht]
\footnotesize
\caption{\footnotesize Order of Convergence of the $L^{\infty}$-error. Results based on experiment 1.}\label{tfive}
\begin{tabular}{|c|c|c|c|} \hline
$C_v$ & HP & normalized EHP with $X=50$ & normalized EHP with $X=6$ \\ \hline
 $C_2$&0.888 &   0.495      &   1.407        \\   
$C_3$&0.812 &     1.160    &   1.381           \\
$C_4$&1.244 &     -0.214   &       -0.498       \\ 
$C_5$& 1.026 &   0.836     &      1.120       \\
$C_6$& 0.957 &    0.396      &    1.166       \\
$C_7$& 1.018 &    0.722   &      0.729        \\
$C_8$& 0.991 &   1.088     &   2.622            \\ \hline 
\end{tabular}
\end{table}

\bibliographystyle{mcom-l}

\end{document}